\documentclass[11pt,a4paper]{article}
\usepackage[linktocpage,colorlinks=true,raiselinks]{hyperref}
\hypersetup{urlcolor=blue, citecolor=red, linkcolor=blue}

\usepackage{comment}
\usepackage{graphicx,fancyhdr}
\usepackage{amsmath,bm}
\usepackage{amssymb}
\usepackage{amsfonts}
\usepackage{amsthm}
\usepackage{xcolor}
\usepackage[T1]{fontenc}
\usepackage{lmodern}
\usepackage{float}

\numberwithin{equation}{section}

\usepackage{paralist}
\numberwithin{equation}{section}

\usepackage{mathtools, mathabx, verbatim}

\def\R{\mathbb{R}}

\def\neweq#1{\begin{equation}\label{#1}}
\def\endeq{\end{equation}}

\newtheorem{thm}{Theorem}[section]
\newtheorem{lem}[thm]{Lemma}
\newtheorem{cor}[thm]{Corollary}
\newtheorem{prop}[thm]{Proposition}

\newtheorem{defn}[thm]{Definition}
\theoremstyle{definition}
\newtheorem{rem}[thm]{Remark}
\theoremstyle{remark}

\newcommand{\ds}{\displaystyle}
\newcommand{\norm}[1]{\left\Vert#1\right\Vert}

\newcommand{\abs}[1]{\left\vert#1\right\vert}

\newcommand{\N}{\mathbb{N}}

\newcommand{\Z}{\mathbb{Z}}

\DeclareMathOperator{\dive}{div}

{\left\{\begin{array}{@{}l@{}}}{\end{array}\right.}

\makeatletter 
\def\@makefnmark{} 
\makeatother 

\title{\bf Energy conservation for 3D Euler and Navier-Stokes
  equations in a bounded domain. Applications to Beltrami
  flows}
\author{Luigi C. Berselli, Elisabetta Chiodaroli, and Rossano Sannipoli
  \\
  Dipartimento di Matematica - Università di Pisa, Italy
  \\
  email: luigi.carlo.berselli@unipi.it, elisabetta.chiodaroli@unipi.it,
  \\
  rossano.sannipoli@dm.unipi.it}
\date{\today}
\begin{document}
\maketitle
\begin{abstract}
  In this paper we consider the incompressible 3D Euler and Navier-Stokes equations in a
  smooth bounded domain. First, we study the 3D Euler equations endowed with slip boundary
  conditions and we prove the same criteria for energy conservation involving the
  gradient, already known for the Navier-Stokes equations. Subsequently, we utilize this
  finding, which is based on a proper approximation of the velocity (and doesn't
 require estimates or additional assumptions on the pressure), to explore energy
  conservation for Beltrami flows. Finally, we explore Beltrami solutions to the
  Navier-Stokes equations and demonstrate that conditions leading to energy conservation
  are significantly distinct from those implying regularity. This remains true even when
  making use of the bootstrap regularity improvement, stemming from the solution being a
  Beltrami vector field.
  \\
  \textbf{MCS:}{ Primary 35Q31; Secondary 76B03.}
\end{abstract}
\section{Introduction}
In this paper we address two problems concerning energy conservation
for weak solutions to incompressible fluids.  Throughout,
$\Omega\subset \mathbb R^3$ will be a bounded domain with strongly
Lipschitz boundary, and $u^E\,, u: (0,T)\times \Omega\to \mathbb R^3$
and $p^{E},\,p: (0,T)\times \Omega\to \mathbb R$ will represent
respectively the velocity vector field of an ideal (or of a viscous)
homogeneous fluids and their associated kinematic pressure.

The first problem pertains to the analysis of energy conservation
for weak solutions to the 3D Euler equations, equipped with a slip
boundary condition at the boundary 
\begin{equation}
\label{eq:Euler}
\begin{cases}
  \partial_t u^E + (u^E\cdot \nabla)\,u^E+\nabla p^{E} = 0 \qquad
  \qquad& \,\text{in}\, (0,T)\times \Omega,
  \\
  \dive u^E = 0 & \,\text{in}\, (0,T)\times \Omega,
  \\
  u^E\cdot n = 0 & \,\text{on}\, [0,T]\times \partial \Omega,
  \\
  u^E(0,x)=u_0^E(x) & \,\text{in}\, \Omega.
\end{cases}
\end{equation}
The second problem focuses on the Leray-Hopf weak solutions to the Navier-Stokes equations (NSE) with Dirichlet boundary condition, that is
\begin{equation}
\label{eq:NSE}
\begin{cases}
  \partial_t u -\Delta u+ (u\cdot \nabla)\,u+\nabla p = 0 \qquad \qquad& \,\text{in}\, (0,T)\times
  \Omega,
  \\
  \dive u = 0 & \,\text{in}\, (0,T)\times \Omega,
  \\
    u = 0 & \,\text{on}\, [0,T]\times \partial \Omega,
    \\
    u(0,x)=u_0(x) & \,\text{in}\, \Omega, 
\end{cases}
\end{equation}
where the kinematic viscosity is set to be equal to 1, without loss of
generality.

It is well known that for smooth solution to~\eqref{eq:Euler} (which
are known to exists only locally in time) the kinetic energy
\begin{equation*}
  E(t) := 
  \frac{1}{2}\|u^E(t)\|_2^2,
\end{equation*}
is constant for $t\in[0,T]$, while for smooth (at least strong)
solutions of the Navier-Stokes equations the following balance
relation holds
\begin{equation}
\label{eq:energy-equality}
    \frac{1}{2}\|u(t)\|_2^2+\int_0^t\|\nabla u(s)\|^2_{2}\,d
    s=\frac{1}{2}\|u_0\|^2_{2}. 
 \end{equation}
Recall also that, for the Euler equations the equality of energy is a
property valid for smooth solutions, which are only known to exist
locally in time, while for weak solutions the balance of energy
cannot even be computed. On the other hand, for (Leray-Hopf) weak
solutions to the Navier-Stokes equations only an inequality in the
balance law is known, which does not exclude the possibility of
having anomalous dissipation.

The investigation of this problem and the verification of the
kinetic energy balance in both scenarios (including also the limit of
vanishing viscosity) have a rich history, tracing back to the works
of Kolmogorov~\cite{Kol1941} and Onsager~\cite{Ons1949}. For a
comprehensive overview, one can refer to
Frisch~\cite{Fri1995}. Subsequently, the interest in this problem
extended to pure mathematicians, particularly after the
contributions of J.L.~Lions and G.~Prodi around 1960. A review of these
developments can be found in~\cite{Ber2021}. Key milestones in
understanding the mathematical subtleties of the problem include
the seminal results by Prodi~\cite{Pro1959} and
Lions~\cite{Lio1960}, as well as that by Constantin, E, and
Titi~\cite{CET1994}, which identified additional conditions on the
velocity field for both the Navier-Stokes equations (NSE) and the
Euler equations.

Over the past 15 years, the interest for the subject has been renewed stimulating several
researchers, who obtained notable improvements.  These include extensions to more general
spaces, as in Cheskidov and Luo~\cite{CL2020}, and explorations of connections with the
singular limit between the viscous and ideal cases, as investigated by Drivas and
Eyink~\cite{DE2019}.

In this paper, our primary focus lies on the implications of assuming conditions on the
gradient of the velocity rather than on the velocity itself. It's worth noting that much
of the existing theory revolves around H\"older, Besov, or fractional spaces of various
kinds, while precise results regarding energy conservation under assumptions on the full
gradient (or on the curl) have only emerged recently.

For the Navier-Stokes equations, precise criteria involving the gradient (as exemplified
in Theorem~\ref{thm:Berselli-Chiodaroli}) have been established in~\cite{BC2020} and by
Beir\~ao and Yang~\cite{BY2019}; related works can also be found in Cheskidov and
Luo~\cite{CL2020}. On the other hand, recent findings regarding the Euler equations
indicate that, more or less, the same conditions applicable to the viscous case are also
valid for the ideal one, as shown in~\cite{BG2024} (for sub-optimal cases) and by Liu,
Wang, and Ye~\cite{LWY2023} (for optimal cases, albeit only for the smallest exponent in
the space variables). The possibility of effectively handling these limit cases became
apparent following the work of De Rosa~\cite{Der2020} and Nguyen, Nguyen, and
Tang~\cite{NgNgT2019}, particularly within functional spaces in which
smooth functions are dense.

Thus, it's expected that in Theorem~\ref{thm:grad-Euler}, we reach the endpoint case (with
further insights can be found in~\cite{BS2024} for fractional spaces). However, the
primary technical advancement here, unlike in other referenced papers, is that we tackle
the problem within a bounded domain. In a bounded domain, conventional techniques
involving smoothing by convolution to approximate the problem are unavailable. Instead, we
employ a modified technique based on transversal mapping to deform the domain, as
elaborated in Section~\ref{sec:preliminaries}. We remark that, prior to our study, very
few results concern the boundary value problem for the Euler equations, see Bardos and
Titi~\cite{BT2018} and Robinson, Rodrigo, and Skipper~\cite{RRS2018}, for the H\"older
case.  Note that our results: a) do not need estimates on the pressure as
in~\cite{BT2018}; b) are not based on reflection and special geometries of the domain as
in~\cite{RRS2018}; c) do not need additional assumptions on the flux as in Drivas and
Nguyen~\cite{DN2018}.

We also recall that many of the results about energy conservation are
intertwined with the Onsager conjecture, as discussed in De Lellis and
Sz\'{e}kelyhidi~\cite{DLS2009}. While the negative aspect of the
Onsager conjecture has been addressed, culminating in the
groundbreaking works by Isett~\cite{Ise2018} and Buckmaster \textit{et
  al.}~\cite{BDLSV2019},  significant efforts still persist
in determining the minimal space-time assumptions necessary for energy
conservation, particularly in the viscous case. Several
recent contributions to this ongoing endeavor include those outlined
in~\cite{BG2024} and in Wang, Wei, Wu, and
Ye~\cite{WWWY2023}.

The first result we prove concerns the conservation of energy under regularity conditions on the gradient of the solution $u^{E}$ to~\eqref{eq:Euler}. 
\begin{thm}
  \label{thm:grad-Euler}
  Let $u_0^E\in H$ and $u^E$ be a weak solution to the Euler
  equations~\eqref{eq:Euler}. If  
  \begin{equation}
    \label{eq:grad-Euler}
    \nabla u^E\in L^\frac{5q}{5q-6}(0,T;L^q(\Omega))\qquad\text{for} \quad
    q\ge\frac{9}{5}, 
\end{equation}
  then, the velocity $u^E$ satisfies the energy equality
\begin{equation*}
    \|u^E(t,x)\|^2_{2}=\|u_0^E(x)\|^2_{2}\qquad\text{for a.e. }\
    t\in[0,T]. 
\end{equation*}
\end{thm}
The result can be also restated in terms of a condition on the  curl
$\omega^E=\nabla\times u^E$.
\begin{cor}
  \label{cor:omega-Euler}
  If $\omega^E\in L^\frac{5q}{5q-6}(0,T;L^q(\Omega))$ for $q\ge\frac{9}{5}$, then energy
  is conserved, provided that the first Betti number of $\Omega$ vanishes.
   \end{cor}
   This Corollary directly follows from Theorem~\ref{thm:grad-Euler}
   by using the estimates from von Wahl~\cite{vWah1992} proving that
   if $\nabla\cdot u^{E}=0$ in $\Omega$ and $u^{E}\cdot n=0$ on
   $\partial \Omega$ then
\begin{equation*}
  \|\nabla u^E\|_p\leq C_p\|\omega^E\|_p\qquad \text{ for }1<p<\infty, 
\end{equation*}
provided that the domain has vanishing first Betti number, i.e.  the
dimension of the homology group $\mathbb{H}^{1}(\Omega,\Z)$ is zero.

Note that the result in Theorem~\ref{thm:grad-Euler} (at least in the
restricted range of $q$ for which is valid) is exactly the same as the one
already known for weak solutions to the Navier-Stokes equations,
cf.~Thm.~\ref{thm:Berselli-Chiodaroli}. (in the latter case observe
that the weak solutions satisfy in addition that the gradient of the
velocity is space-time
square-integrable). Theorem~\ref{thm:grad-Euler} extends the results
in Liu, Wang, and Ye~\cite{LWY2023} where only the limit case $q=9/5$
was considered. Moreover and more significantly, we can handle not
only the Cauchy problem, but also the boundary value problem in a
smooth domain $\Omega\subset\R^3$, differently from similar results
from Wang~\textit{et al.}~\cite{WWWY2023}.
\medskip

In the second part of the paper we further investigate the implications of
Theorem~\ref{thm:grad-Euler} on the energy conservation of Beltrami
(also known as Trkal) flows, which are families of solutions
characterized by a specific geometric constraint. Beltrami solutions
play a significant role in fluid dynamics, as they represent a set of stationary
solutions to the Euler equations~\eqref{eq:Euler}. Specifically, these
solutions are characterized by the property that $\omega^E$ is
proportional to the velocity field itself:
\begin{equation}
  \label{eq:Beltrami}
    \omega^E(x,t) = \lambda(x,t)u^E(x,t),
\end{equation}
where $\lambda(\cdot,\cdot)$ is a suitable scalar function of the time and/or space
variables.

The simplest and smoothest case corresponds to $\lambda \equiv 0$, i.e. potential flows are a particular case of Beltrami flows. We remark that Beltrami flows are genuinely 3D flows
 since in the 2-dimensional setting $\omega$ is orthogonal to the plane of motion. Additionally, by
  employing the \textit{Lamb vector} $\omega^E\times u^E$, it is possible to express the convective term in the
  following \textit{rotational formulation} 
\begin{equation*}
  (u^E\cdot\nabla)\, u^E=\omega\times u^E+\frac{1}{2}\nabla|u^E|^2.
\end{equation*}
Considering this formulation for Beltrami flows, the convective term amounts to a gradient
(the Bernoulli pressure), which can be included into the pressure. Formally, Beltrami
flows satisfy linear (non-local) transport or Stokes equations, implying that the flow is laminar.
However there are two warnings: i) the numerical treatment of the pressure, and especially
that of the Bernoulli pressure, leads to a stiff problem. If not using \textit{pressure-robust}
numerical methods, the results at very high Reynolds numbers could be affected by
instabilities (see Gauger, Linke, and Schroeder~\cite{GLS2019}); b) from a formal point of view  $(\omega^E\times u^E)\cdot u^E=0$ almost everywhere in $(0,T)\times\Omega$; however, for weak
solutions $\omega^E$ is only a distribution and  this equality is not enough to imply that the  integral
$\int_0^T\int_{\Omega}(\omega^E\times u^E)\cdot u^E\, dx dt$ is well-defined. See a similar discussion
in~\cite{BS2024}.

Extending also to the viscous case results already sketched in~\cite{BS2024}, we present
some elementary observations on the regularity implied by the geometric
constraint~\eqref{eq:Beltrami}. If $\lambda(x,t)\equiv \lambda \in \mathbb{R}$ (a
circularly polarized plane-wave), then $u^{E}$ is smooth, and the
conservation of energy follows. This can be established through a
standard bootstrap argument, assuming regularity conditions on the boundary of $\Omega$
and $\mathbb{H}^{1}(\Omega,\Z)=0$. Consequently, a continuation argument for smooth
solutions holds, provided that the initial datum is smooth

The second observation stems from a straightforward computation when
$\lambda(x,t)=\lambda(t)\in L^{p}(0,T)$, for some $p\geq 1$. Note that in this scenario,
$\nabla\cdot\omega^E=\lambda(t)(\nabla\cdot u^E)=0$, thereby satisfying the divergence-free
constraint without requiring any additional assumption on $\lambda(t)$. Consequently,
from~\eqref{eq:Beltrami}, we deduce that $\omega^E \in L^p(0,T;L^2(\Omega))$, thereby
implying better regularity on $u^E$. Specifically for the Euler equations, this implies
$u^E\in L^p(0,T;H^1(\Omega))$. By iterating this procedure, we establish that if
$\lambda(t)\in L^p(0,T)$ for some $p\geq 1$, and $u^E\in L^\infty(0,T;L^2(\Omega))$,
then
\begin{equation*}
    \omega^E\in L^\frac{p}{3}(0,T;H^2(\Omega))\xhookrightarrow{}L^\frac{p}{3}(0,T;L^\infty(\Omega)).
\end{equation*}
Hence, if $p\geq3$,  the Beale-Kato-Majda~\cite{BKM1984} criterion for
continuation of smooth solutions holds.

In the case of the Navier-Stokes equation, from $\omega \in L^p(0,T;L^2(\Omega))$, we
deduce
\begin{equation*}
u\in L^p(0,T;H^1(\Omega))\hookrightarrow L^p(0,T;L^6(\Omega)),
\end{equation*}
thus, if $p\geq 2$, the next iteration gives $\omega \in L^{p/2}(0,T;L^6(\Omega))$ and
consequently $\nabla u\in L^{p/2}(0,T;L^6(\Omega))$. By applying the scaling invariance
criterion: if
\begin{equation}
\label{eq:scaling-nabla}
\nabla u\in L^{r}(0,T;L^{s}(\Omega))\quad\text{for }
\frac{2}{r}+\frac{3}{s}=2,\text{ with } s>\frac{3}{2},
\end{equation}
then $u$ is a strong solution to the NSE (cf.~Beir\~ao da Veiga~\cite{Bei1995a}
and~\cite{Ber2002a} for the problem in a bounded domain), we obtain that if $p\geq\frac{8}{3}$, then 
\begin{equation*}
\frac{2}{p/2}+\frac{3}{6}\leq 2.
\end{equation*}

Therefore, an initial elementary result is as follows.
\begin{prop}
  \label{prop:Elementary}
  Let $u^E$ be a weak solution to the Euler equations~\eqref{eq:Euler}, which is a
  Beltrami field with $\lambda\in L^p(0,T)$, for $p\geq3$. Let
  $u_0^E\in H^3(\Omega)\cap V_\tau$, then $u^{E}$ is the unique classical solution
  of~\eqref{eq:Euler} in $[0,T]$ and satisfies the energy equality.

  Let $u$ be a weak solution to the Navier-Stokes equations~\eqref{eq:NSE}, which is a
  Beltrami field with $\lambda\in L^p(0,T)$, for $p\geq8/3$. Let $u_0\in V_0$, then $u$ is
  the unique strong solution of~\eqref{eq:NSE} in $[0,T]$ and satisfies the energy
  equality.
\end{prop}
These two elementary examples also illustrate how it would be straightforward to fall into
a regularity class even under mild hypotheses on $\lambda$. Our goal is to identify
classes in which the energy is conserved without being trivially categorized into a
smoothness class. 

In cases where $\lambda$ depends on spatial variables, maintaining the divergence
condition requires $\nabla\lambda\cdot u^{E}=\nabla\lambda\cdot u^E=0$. This condition
influences the effective velocity fields, especially in classical solutions, as discussed
by Beltrami~\cite{Bel1873}. Recent research by Enciso and Peralta~\cite{EP2016} and
Abe~\cite{Abe2022} explores the potential existence of non-trivial Beltrami
fields. Hereafter, we consider weak solutions $u^{E}$ and $u$ which are Beltrami fields. Concerning Beltrami solutions to the Euler equations in three-dimensional torus we refer also to~\cite{BS2024}, where sufficient conditions in fractional Sobolev spaces are found on $\lambda$ in order to have conservation of energy.
\medskip

The first result of this paper, regarding the energy equality for solutions to problem \eqref{eq:Euler} which are Beltrami fields as well, is the following. 
\begin{thm}
\label{thm:BeltramifieldthmEuler}
Let $u^E$ be a weak solution to the Euler equation~\eqref{eq:Euler},
which is a Beltrami field, i.e. \eqref{eq:Beltrami} is satisfied.  Let
$\lambda \in L^\alpha(0,T;L^\beta(\Omega))$, where
\begin{equation*}
  \beta =\begin{cases}
       \frac{6\alpha}{2\alpha-5} & \text{ if }\alpha >
       \frac{5}{2},
       \\
      \infty &\text{ if } \alpha = \frac{5}{2}.
  \end{cases} 
\end{equation*}
Then, the velocity $u^{E}$ satisfies the energy equality for a.e. $t\in[0,T]$.
\end{thm}
\bigskip

In the last part of this paper we study the energy conservation
results for Leray-Hopf weak solutions of the Navier-Stokes
equations~\eqref{eq:NSE}, when they are Beltrami fields, as well.

We recall that regarding the Navier-Stokes equations, criteria for energy conservation
have been established in~\cite{BC2020} and subsequently refined and extended to
non-Newtonian fluids by Beir\~ao and Yang~\cite{BY2019}.
\begin{thm}
  \label{thm:Berselli-Chiodaroli}
Let $u$ be a Leray-Hopf weak solution
  of~\eqref{eq:NSE}.  Let us
  assume that one the following conditions is satisfied
  \begin{itemize}
  \item[(i)]  $\nabla u\in
  L^\frac{q}{2q-3}(0,T;L^q(\Omega))$,\qquad for\quad  $\frac{3}{2}< q<\frac{9}{5}$;
  \item[(ii)] $\nabla u\in L^\frac{5q}{5q-6}(0,T;L^q(\Omega))$,\qquad for
    \quad$q\ge\frac{9}{5}$, 
 \end{itemize}
 Then, the velocity $u$ satisfies the energy equality~\eqref{eq:energy-equality}.
\end{thm}
\begin{rem}
  Also in this case one can state the result using $\omega$ instead of $\nabla u$, as in
  Corollary~\ref{cor:omega-Euler}. Since the full trace of velocity is zero at the
  boundary the geometric conditions can be relaxed and a smooth bounded domain is enough.
\end{rem}
Starting from Theorem~\ref{thm:Berselli-Chiodaroli} we provide the counterpart of
Theorem~\ref{thm:BeltramifieldthmEuler} for the Navier-Stokes equations.
\begin{thm}
\label{thm:Beltrami-NSE}
    Let $u$ be a Leray-Hopf weak solution of~\eqref{eq:NSE}, which is
    a Beltrami field. Let 
    $\lambda \in L^\alpha(0,T;L^\beta(\Omega))$, where
\begin{equation*}
  \beta =\begin{cases}
       \frac{6\alpha}{2\alpha-5} & \text{ if } \alpha >
       \frac{5}{2},
       \\
      \infty &\text{ if } \alpha = \frac{5}{2}.
  \end{cases}  
\end{equation*}
Then, the velocity $u$ satisfies the energy equality for a.e. $t\in[0,T]$.
\end{thm}
\begin{rem}
  The conditions on $\lambda$ remain the same as in the case of the Euler equations, as
  the bootstrap process starts from $u\in L^\infty(0,T;L^2(\Omega))$ in both scenarios,
  and utilizing the information $u\in L^2(0,T;L^6(\Omega))$ doesn't appear to enhance the
  regularity. What is relatively different in this case is that for the Euler equations,
  achieving $\omega\in L^1(0,T;L^\infty(\Omega))$ ensures regularity. However, this
  is unattainable unless $\lambda$ itself is bounded in the spatial variables. In
  contrast, for the Navier-Stokes equations, much less restrictive assumptions on
  $\lambda$ suffice to fall into classes of uniqueness and regularity as described
  in~\eqref{eq:scaling-nabla}.
\end{rem}
The arguments that lead to Theorem~\ref{thm:Beltrami-NSE}, when combined with classical
scaling invariant conditions for regularity of weak solutions, allow us to prove also the
following result on the regularity of Beltrami fields. Before giving its precise statement
we need to introduce first the following decomposition of the half-line
$(3,+\infty)=\bigcup_{n\in \mathbb N} I_n$, where $I_n= L_n \cup R_n$, and
\begin{equation}
  \label{eq:LnRn}
    L_n := \begin{cases}
      (12,24] & n=1
      \\
      \\
       \left(\frac{6(n+1)}{2n-1},\frac{6(n+1)^2}{2n^2+n-2}\right] & n\ge 2,
    \end{cases}
    \qquad
    R_n := \begin{cases}
      (24,+\infty) & n=1
      \\
      \\
       \left(\frac{6(n+1)^2}{2n^2+n-2},\frac{6n}{2n-3}\right] & n\ge 2.
    \end{cases}
\end{equation}
\begin{thm}
  \label{thm:regularityNSE}
  Let $u$ be a Leray-Hopf weak solution of~\eqref{eq:NSE} corresponding to $u_0\in V_{0}$
  as initial datum and let $u$ be also a Beltrami field.  For any fixed $\beta>3$, there
  exists a unique $\overline{n}=\overline{n}(\beta)\in \mathbb{N}$, such that if $\lambda \in
  L^{\alpha}(0,T;L^\beta(\Omega))$ and 
$\alpha$ satisfies
\begin{equation}
  \label{eq:regularityLnRn}
    \begin{cases}
        
    \displaystyle\frac{2}{\alpha}+\frac{3}{\beta} =
    1-\frac{1}{2(\overline{n}+2)} & \quad\text{for}\,\,\, \beta\in
    L_{\overline{n}},
    \\
    \\
    \displaystyle\frac{2}{\alpha}+\frac{3}{\beta} =
    1-\frac{1}{2(\overline{n}+1)} &\quad \text{for}\,\,\, \beta \in R_{\overline{n}},
    \end{cases}
\end{equation}
then $u$ is the unique strong solution in $(0,T)$.
\end{thm}
As a corollary of the above result we have the following scaling invariant criterion in
terms of $\lambda$, for regularity of weak solutions.
\begin{cor}
\label{cor:regularityNSE}
Let $u$ be a Leray-Hopf weak solution of~\eqref{eq:NSE} corresponding to $u_0\in V_{0}$ as
initial datum and let $u$ be also a Beltrami field. If
$\lambda \in L^{\alpha}(0,T;L^\beta(\Omega))$, where $\alpha$ and $\beta$ satisfy
  \begin{equation*}
      \frac{2}{\alpha}+\frac{3}{\beta }<1,
\end{equation*}
then there exists $\beta_0$ such that, for every
$\beta \in (3,\beta_0)$, $u$ is the unique strong solution in $(0,T)$.
\end{cor}

\bigskip

\textbf{Plan of the paper:} In Section~\ref{sec:preliminaries} we give
the definitions of weak solutions to equations
\eqref{eq:Euler}-\eqref{eq:NSE}, introduce the functional spaces we
will work with, and discuss about the mollification procedures in
space and time, discussing the connected results that will be helpful
in the rest of the paper. In Section~\ref{sec:weak} we will give the
proofs of Theorems
\ref{thm:grad-Euler}-\ref{thm:BeltramifieldthmEuler} and in
Section~\ref{sec:BeltramiNSE} we prove Theorem
\ref{thm:Beltrami-NSE}-\ref{thm:regularityNSE}.

%

\section{Functional setting, weak solutions and mollifiers}
\label{sec:preliminaries}
In this section we first introduce the notation and the precise
definitions of solutions we want to deal with. Then, we compare our
results with the ones in the existing literature.  We will use the
customary Sobolev spaces $(W^{k,p}(\Omega),\,\|\,.\,\|_{W^{k,p}})$ and
we denote the $L^p$-norm by $\|\,.\,\|_p$. %
We will not distinguish scalar and vector valued spaces,
since this will be clear from the context.

For  a Banach space $X$ we will also denote the
usual Bochner spaces of functions defined on $[0,T]$ and with values in $X$ by
$(L^p(0,T;X),\|\,.\,\|_{L^p(X)})$. In the case $X=L^q(\Omega)$ we denote the norm of
$L^{p}(0,T;L^{q}(\Omega))$ simply by $\|\,.\,\|_{p,q}$.
\subsection{On the weak solutions of the Euler and Navier--Stokes equations}
\label{sec:weak-solutions}
For the weak formulation of 3D Euler equations~\eqref{eq:Euler} and Navier--Stokes equations~\eqref{eq:NSE}, we introduce the following spaces 
\begin{equation*}
\begin{aligned}
  &H_\tau= \{v \in L^2(\Omega): \nabla\cdot v =0\;\text{in}\;
  \Omega,\,v\cdot n =0 \;\text{on}\;\partial\Omega\},
  \\
  &V_\tau= \{v \in H^1(\Omega): \nabla\cdot v =0\;\text{in}\;
  \Omega,\,v\cdot n =0 \;\text{on}\;\partial\Omega\}.
\end{aligned}
\end{equation*}
The Hilbert space $H_\tau$ is endowed with
the natural $L^2$-norm $\norm{\,.\,}_{2}$ and inner product $(\cdot,\cdot)$, while $V_\tau$ with
the norm $\norm{\nabla v}_{2}$ and inner product $((u,v)):=(\nabla u,
\nabla v)$.

The space of test functions considered to define a weak solution of
the Euler equation~\eqref{eq:Euler} is the following
\begin{equation*}
  \mathcal{D}_T = \{\varphi \in C_{0}^{\infty}([0,T[; C^\infty(\overline{\Omega})) :
  \nabla\cdot\varphi =0\;\text{in}\; \Omega,\,\varphi\cdot n =0 \;\text{on}\;\partial\Omega\}. 
\end{equation*}
\begin{defn}[Euler weak solutions]
\label{def:weak-E}
  Let $u_0^E\in H_\tau$. A measurable function $u^E:(0,T)\times \Omega\to \mathbb R^3$ is called a
  weak solution to the Euler equation if $u^E \in C(0,T;w\text{-} H_\tau)$ (continuous w.r.t the weak topology) is such that
  \begin{equation*}
    \int_0^{T} [(u^E,\partial_t\varphi)+((u^E\otimes u^E),\nabla \varphi)]\,dt =
    -(u^E_0,\varphi(0)) \qquad  \forall\, \varphi \in \mathcal{D}_T. 
    \end{equation*}
\end{defn}
For the functional setting of the Navier-Stokes
equations~\eqref{eq:NSE} we introduce also the space $H$ as the closure
of the space $\mathcal{V}$ of smooth and divergence-free vectors
fields, with compact support in $\Omega$ with respect to the norm of
$L^2(\Omega)$ as well as the space $V_{0}$ as the closure
of the space $\mathcal{V}$ of smooth and divergence-free vectors
fields, with compact support in $\Omega$ with respect to the norm of
$W^{1,2}_{0}(\Omega)$.

\begin{defn}[Leray-Hopf  weak solutions]
\label{def:weak-NS}
  A vector field $u\in L^\infty(0,T; H)\cap L^2(0,T;V_0)$ is a Leray-Hopf weak solution to
  the Navier-Stokes equations~\eqref{eq:NSE} if
\begin{itemize}
	\item [(i)]  $u$ is a solution of ~\eqref{eq:NSE} in the weak sense, i.e.
  \begin{equation*}
    \int_0^T(u,\partial_t \phi)-(\nabla u,\nabla \phi)-((u\cdot\nabla)\,
    u,\phi)\,\,dt=-(u_0,\phi(0)), 
  \end{equation*}
  for all $\phi\in C^\infty_0([0,T[\times \Omega)$ with $\nabla\cdot \phi=0$;
\item [(ii)] $u$ satisfies the global energy inequality
  \begin{equation*}
    \frac{1}{2}\|u(t)\|^2_{2}+\int_0^t\norm{\nabla
      u(s)}^2_{2}\, ds\leq\frac{1}{2}\|u_0\|^2_{2}\qquad 
    \forall\,t\geq0;
  \end{equation*}
\item [(iii)]  the initial datum is attained in the strong sense of $L^2(\Omega)$
  \begin{equation*}
    \|u(t)-u_0\|\to0\qquad\text{as}\quad  t\to0^+.
  \end{equation*}
\end{itemize}
\end{defn}
\subsection{Mollification in space and time}\label{subsec:molspacetime}
As usual, the proof of Theorem~\ref{thm:grad-Euler} is mainly based on
estimates on the non-linear term in the energy balance, by means of
suitable mollification in space and time, separately. This is done in
order to make the formal multiplication by $u^E$ and integration by
parts rigorously justified. To this, end we need a proper way to smooth
a divergence-free vector field preserving both the boundary values and
the incompressibility.
\subsubsection{Mollification in space}\label{molspace} 
For $V_{0}$ a standard density argument is enough to reach the desired
approximation. On the other hand, for the Euler equations, the
condition $u\cdot n=0$ on the boundary makes the approximation a bit
more complicated. That's why we need to find another way to
approximate a weak solution to the Euler equations with smooth
functions that preserve both the zero divergence condition and the (at
least the slip) boundary condition. This kind of approximation is
given in~\cite{EG2016} where the authors construct suitable mollifiers
that preserve boundary conditions. Here we briefly give the main
definitions and results and we refer to~\cite{EG2016} (and the
references therein) for the details.

Let $\Omega\subset \mathbb R^3$ be an open, bounded and strongly
Lipschitz connected set. Since $\Omega$ is bounded, we can find
$x_0\in \Omega$ and $r_0>0$, such that $\Omega\subset B_{r_0}(x_0)$,
where $B_{r_0}(x_0)$ is a ball centered in $x_0$ with radius
$r_0>0$. Then, the open set
$\mathcal{O}= B_{r_0}(x_0)\setminus \overline{\Omega}$ is a strongly
Lipschitz domain and it is possible to prove the existence of a vector
field $V\in C^\infty(\mathbb R^3)$ which is globally transversal for
$\mathcal{O}$, i.e. $\exists \,\alpha >0$, such that
$n(x)\cdot V(x)\ge \alpha$ for a.e. $x\in \partial\mathcal{O}$, and
its euclidean norm satisfies $\|V(x)\|_{l^2}=1$, for all
$x\in \partial \mathcal{O}$. We then define the map
\begin{equation*}
    \theta_{\delta} (x) : x\in \mathbb R^3 \to x+\delta V(x)\in \mathbb R^3,
\end{equation*}
with the following properties:
$\theta_{\delta}\in C^\infty(\mathbb{R}^{3})$ for all
$\delta \in [0,1]$ and there exists a $\xi>0$ such that
\begin{equation*}
    \theta_{\delta}(\mathcal{O}) +B_{3\delta\xi }(0)\subset
    \mathcal{O}\, \qquad \forall\, \delta \in [0,1], 
\end{equation*}
where $B_{3\delta\xi}(0)$ is a ball centered at the origin of radius
$3\delta\xi$.
\\

Let $v \in L^1(\Omega,\mathbb R^3)$ and let $\overline{v}$ be the
extension to zero over off $\mathcal{O}$. Now we can define the
desired mollifier as follows
\begin{equation*}
(\mathcal{K}^{div}_\delta v)(x) :=
\int_{B}\rho(y)\det(J_\delta(x))J^{-1}_\delta(x)\,\overline{v}(\theta_{\delta}(x)+(\delta\,\xi)y)\,dy, 
\end{equation*}
where $B$ is the unit ball centered at the origin, $\rho(y)$ is the
standard Friedrichs  mollifier and $J_\delta$ is the Jacobian of the
map $\theta_{\delta}$. Observe that $J_\delta$  satisfies for all $m\in
\mathbb N$ and some positive constant $c>0$ 
\begin{equation}
  \label{supjacobian}
    \sup_{x\in \Omega} \|D^m(\det(J_\delta(x))J^{-1}_\delta(x)-\mathbb I)\|_{l^2}\le c\,\delta,
\end{equation}
where $\mathbb I\in \mathbb R^{3\times 3}$ is the identity matrix in $\mathbb R^3$. 
In particular we have
\begin{lem}
  \label{lem:smoothness}
For all $\delta\in(0,1]$ it holds that $\mathcal{K}^{div}_\delta :
L^1(\Omega;\mathbb R^3)\to C^\infty_0(\Omega;\mathbb R^3)$ and
$\mathcal K^\text{div}_\delta v $ admits a continuous extension to $\overline{\Omega}$.
\end{lem}
Defining for $1\le q \le \infty$ 
\begin{equation*}
\begin{aligned}
V_{div,q}(\Omega) = \{v \in L^q(\Omega), \nabla\cdot v \in
L^q(\Omega)\},
\\
\widetilde{V}_{div,q}(\Omega) = \{v \in L^q(\Omega), \nabla\cdot \overline{v} \in L^q(\mathbb R^3)\},
\end{aligned}
\end{equation*}
the fundamental property of this operator is that $\nabla\cdot
\mathcal{K}^{div}_\delta v=0$ if $\nabla\cdot v=0$ and $v\cdot n=0$ 
( hence $\overline{v}\in \widetilde{V}_{div,q}$) and the following general convergence holds true.
%
\begin{lem}\label{lem:convergenceq}
  There exists a $\delta_0$, such that the sequence
  $\mathcal{K}^{div}_\delta$, with $\delta \in [0,\delta_0]$, is
  uniformly bounded in $\mathcal{L}(L^q(\Omega),L^q(\Omega))$, for all
  $1\le q<\infty$. Moreover, we have
  \begin{equation*}
    \label{eq:convergenceLq}
        \lim_{\delta \to 0}\|\mathcal{K}^{div}_\delta
        v-v\|_{q}=0\, \qquad \forall\,v \in
        L^q(\Omega), 
    \end{equation*}
    and
    \begin{equation*}
        \lim_{\delta \to 0}\|\nabla \cdot(\mathcal{K}^{div}_\delta
        v-v)\|_{q}=0\, \qquad \forall\,v \in
        \widetilde{V}_{div,q}(\Omega). 
    \end{equation*}
\end{lem}


For the purposes of this paper, we need some boundedness result
concerning the (full) gradient of $\mathcal{K}^{div}_\delta v$, which
is not contained in~\cite{EG2016}. In particular we prove the
following result.
\begin{lem}
  \label{lem:convergeneW1q}
  For every $v\in W^{1,q}(\Omega)$ and $\delta \in (0,1]$, there
  exists a positive constant $C>0$ such that
    \begin{equation*}
        \| \nabla\mathcal{K}^{div}_\delta v\|_{q}\le C.
    \end{equation*}
    
\end{lem}
\begin{proof}
   Let us first compute the gradient of $K_\delta^{div}v$
   \begin{equation*}
   \begin{split}
       \nabla K_\delta^{div}v =
       &\nabla(\det(J_\delta(x))J^{-1}_\delta(x))\int_{B}\rho(y)\,\overline{v}(\theta_{\delta}(x)+(\delta
       r)y)\,dy \\ 
       &+\det(J_\delta(x))\int_{B}\rho(y)\,\nabla\overline{v}(\theta_{\delta}(x)+(\delta r)y)\,dy.
       \end{split}
   \end{equation*}
Considering the $L^q$-norm, using Minkowski inequality, and since
$v\equiv \overline{v}$ in $\Omega$, we have 
\begin{equation*}
    \begin{split}
       \|\nabla K_\delta^{div}v\|_{q} \le & \bigg(
       \int_\Omega
       \Big\vert\nabla(\det(J_\delta(x))J^{-1}_\delta(x))\int_{B}\rho(y)v(\theta_{\delta}(x)+(\delta
       r)y)\,dy\Big\vert^q\,dx\bigg)^{\frac{1}{q}}
       \\
       &+\bigg(\int_\Omega\Big\vert\det(J_\delta(x))\int_{B}\rho(y)\nabla
       v(\theta_{\delta}(x)+(\delta
       r)y)\,dy\Big\vert^q\,dx\bigg)^\frac{1}{q}. 
       \end{split}
\end{equation*}
By virtue of~\eqref{supjacobian},  for every $\delta \in (0,1]$ we get
\begin{equation*}
     \sup_{x\in \Omega} \abs{\nabla(\det(J_\delta(x))J^{-1}_\delta(x))}\le c\,\delta \le c.
   \end{equation*}
 Moreover, by definition of $\theta_{\delta}$, for every $\delta \in (0,1]$, the following estimate holds
 \begin{equation*}
     \sup_{x\in \Omega} \abs{\det(J_\delta(x))}\le 1+c\,\delta \le 1+c.
   \end{equation*}
   Eventually, by Jensen's inequality and the properties of the
   Friedrichs  mollifier we can write
\begin{equation*}
 \|\nabla K_\delta^{div}v\|_{q} \le c\|v\|_{q}+(1+c)\|\nabla v\|_{q} \le C ,  
\end{equation*}
for every $\delta \in (0,1]$, ending the proof.
\end{proof}
Obviously triangular inequality implies the following result
\begin{cor}\label{cor:boundednessgrad}
    Let $v \in W^{1,q}(\Omega)$, then there exists a positive constant $\Bar{C}>0$, such that
    \begin{equation*}
        \lim_{\delta \to 0} \|\nabla K_\delta^{div}v-\nabla v \|_{q}\le \Bar{C}.
    \end{equation*}
\end{cor}
%
\begin{rem}
  If $v\in L^q(\Omega)$, Lemma~\ref{lem:convergenceq} guarantees the
  $L^q$-convergence of the regularized function $K^{div}_\delta v$ to
  $v$. Let us stress that if $v\in W^{1,q}(\Omega)$, the
  $W^{1,q}$-convergence result does not hold in general. Indeed, let
  us suppose that $v\in W^{1,q}(\Omega)\cap H$ and let us suppose that
  $K^{div}_\delta v\to v$ in $W^{1,q}(\Omega)$. Since
  $K^{div}_\delta v\in C^{\infty}_0(\Omega)$, then by trace inequality
  we would have
  \begin{equation*}
    \|v\|_{L^q(\partial\Omega)}=\|K^{div}_\delta
    v-v\|_{L^q(\partial\Omega)}\le C\|
    K^{div}_\delta v-v\|_{W^{1,q}(\Omega)}\overset{\delta\to0}{\to} 0,
    \end{equation*}
    implying $v=0$ almost everywhere on $\partial \Omega$. But this is
    not necessarily true since $v$ satisfies $v\cdot n =0$ on
    $\partial \Omega$.\\

    The bounds proved show that
    $\nabla K^{div}_\delta v \rightharpoonup \nabla v$ weakly in
    $L^{q}(\Omega)$, but the latter argument, in particular, shows that
    in general strong convergence does not holds. One could prove
    strong convergence in $W^{s,p}(\Omega)$ for all $0\leq s<1/p$,
    that is for fractional spaces on which the trace operator is not
    defined.
  \end{rem}
The last result we prove in this subsection concerns the uniform
continuity in time of $K^{div}_\delta v$ when $v$ is a function
defined on $[0,T]\times \Omega$, and smoothing concerns only the space
variables. 
\begin{prop}
  \label{Prop:uniformconvergenceintime}
    Let $v \in C([0,T],L^q(\Omega))$, then for every $q\in
    [1,+\infty)$ it holds 
    \begin{equation*}
        \lim_{\delta\to 0} \|K_\delta^{div}v-v\|_{\infty,q}=0.
    \end{equation*}
\end{prop}
\begin{proof}
    For all fixed $t\in [0,T]$, by Lemmata~\ref{lem:smoothness}
    and~\ref{lem:convergenceq} it follows that 
    \begin{equation*}
      K_\delta^{div}v(x,t) \in C^\infty_0(\Omega), \qquad \text{and}
      \qquad K_\delta^{div}v(x,t) \stackrel{\delta\to 0}{\to}v(x,t)
       \text{ in }L^q(\Omega), 
    \end{equation*}
    and we want to prove that the convergence is also uniform in time.
    By~\eqref{supjacobian}, it follows that there exists a positive constant $C>1$
    such that
    \begin{equation*}
\sup_{x\in\Omega}\|\det (J_\delta(x))J^{-1}_\delta(x)-\mathbb I\|_{l^{2}} \le C.
    \end{equation*}
    From the properties of Friedrichs  mollifier and the definition
    of the operator $K_\delta^{div}$, we have for every $q\ge 1$ and for all
    fixed $t \in [0,T]$
\begin{equation}\label{eq:estimatekdeltadivq}
\|K_\delta^{div}v(t)\|_{q} \le C \|v(t)\|_{q}.
\end{equation}
By the uniform continuity in $[0,T]$ of $v$ it follows that $q \in [1,+\infty)$ 
\begin{equation}\label{eq:uniformcontinuity}
    \forall\,\lambda >0, \,\exists\,\delta>0:\ |t-s|<\delta \implies \|u(t)-u(s)\|_{q}<\frac{\lambda}{3C}. 
\end{equation} 
Hence, let us choose a partition of the interval $[0,T]$, let us say
$0=t_0<t_1<\cdots<t_N= T$, such that $|t_{i+1}-t_i|< \delta$, for all
$i=0,...,N-1$. Then, for any $t \in [0,T]$, there exists an index
$i_0\in \{0,...,N\}$, such that
\begin{equation*}
    |t-t_{i_0}|< \delta.
\end{equation*}
Therefore by triangular
inequality,~\eqref{eq:estimatekdeltadivq}-\eqref{eq:uniformcontinuity},
and the fact that $C>1$, we get
\begin{equation*}
\begin{split}
\|K_\delta^{div}v(t)-&v(t)\|_{q} \le
\|K_\delta^{div}\big(v(t)-v(t_{i_0})\big)\|_{q}
\\
&+ \|K_\delta^{div}v(t_{i_0})-v(t_{i_0})\|_{q} +
\|v(t_{i_0})-v(t)\|_{q}
\\
&\le
(C+1)\|v(t_{i_0})-v(t)\|_{q}+\|K_\delta^{div}v(t_{i_0})-v(t_{i_0})\|_{q}
\\
&< \frac{C+1}{3C}\lambda + \|K_\delta^{div}v(t_{i_0})-v(t_{i_0})\|_{q}
\\
&<\frac{2}{3}\lambda + \|K_\delta^{div}v(t_{i_0})-v(t_{i_0})\|_{q}.
\end{split}
\end{equation*}
Now, if we fix $\delta_0>0$ such that
\begin{equation*}
    \max_{i=0,...,N}\|K_\delta^{div}v(t_{i})-v(t_{i})\|_{q}<
    \frac{\lambda}{3}, \qquad \forall\,\delta \in (0,\delta_0), 
\end{equation*}
then we finally get
\begin{equation*}
\|K_\delta^{div}v(t)-v(t)\|_{q}<\lambda,
\end{equation*}
and the conclusion follows.
\end{proof}

\subsubsection{Mollification in time}\label{moltime}
Mollification in time is a standard Friedrichs mollification. Indeed
we consider the (time) mollification operator, denoted in the sequel
by $(\cdot)_\varepsilon$, defined for a space-time function
$\Phi:\,(0,T)\times \Omega\to\R^3$ by
\begin{equation*}
  (\Phi)_{\varepsilon}(t,x) :=\int_0^{T}
  k_\varepsilon(t-\tau)\Phi(\tau,x)\,d\tau\qquad \text{for }0<\varepsilon<T, 
\end{equation*}
where $k$ is a $C^\infty_0(\R)$, real-valued, non-negative
even function, supported in $[-1, 1]$, with $\int_{\R} k(s)\,\ds=1$,
and $k_\varepsilon(t):=\varepsilon^{-1}k(t/\varepsilon)$ (standard Friedrichs mollification with respect to the time
variable).

We end this section considering some results that connect the two
mollification operators.
\begin{lem}\label{lem:commutation}
  Let $v \in L^1(0,T; L^1(\Omega))$, then the mollification operators
  in space and time commute, i.e.
    \begin{equation*}
        (K_\delta^{div}v)_\varepsilon= K_\delta^{div}(v_\varepsilon).
    \end{equation*}
\end{lem}
\begin{proof}
    It is just an application of Fubini's Theorem, indeed
    \begin{equation*}
    \begin{split}
(K_\delta^{div}v)_\varepsilon&= \int_0^T k_\varepsilon(t-\tau)\int_{B}\rho(y)\det(J_\delta(x))J^{-1}_\delta(x)\overline{v}(\theta_{\delta}(x)+(\delta r)y,\tau)\,dy\,d\tau\\
&=\int_B \rho(y) \det(J_\delta(x))J^{-1}_\delta(x)\int_0^T k_\varepsilon(t-\tau)\overline{v}(\theta_{\delta}(x)+(\delta r)y,\tau)\,d\tau\,dy \\
&= K_\delta^{div}(v_\varepsilon).
\end{split} \end{equation*} 
\end{proof}
As a consequence of Lemma~\ref{lem:commutation},
Corollary~\ref{cor:boundednessgrad}, and
Proposition~\ref{Prop:uniformconvergenceintime}, we have the following
two results.
\begin{cor}\label{cor:boundednessgradtimespace}
    Let $v \in L^p(0,T; W^{1,q}(\Omega))$ for any $p,q \in (1,+\infty)$. Then there exists a positive constant $\Tilde{C}>0$, such that for any fixed $\varepsilon>0$, we have
\begin{equation*}
    \|\nabla (K_\delta^{div}v)_\varepsilon-\nabla(v)_\varepsilon\|_{p,q}\le\Tilde{C}.
\end{equation*}
\end{cor}
\begin{cor}\label{cor:LinftyLq}
Let $v \in L^p(0,T; L^q(\Omega))$ for any $p,q \in (1,+\infty)$. Then for any fixed $\varepsilon>0$ we have
\begin{equation*}
    \lim_{\delta\to0^+}\|(K_\delta^{div}v)_\varepsilon-(v)_\varepsilon\|_{\infty,q}=0.
\end{equation*}
\end{cor}
\section{On the energy equality for Euler weak solutions.}\label{sec:weak}
Before starting the proof we recall the following property of weak
solutions.  Let $u^E_0\in H$ and let $u^E$ be a weak solution
of~\eqref{eq:Euler}. Then, since $u^E$ is weakly continuous in time, we
consider it as being already redefined on a set of zero Lebesgue
measure in such a way that $u^E(t)\in L^2(\Omega)$ for all $t\in[0,T)$
and it satisfies the identity
  \begin{equation}
  \label{eq:Eulerequality}
  (u^E(s), \varphi(s))= (u^E_0, \varphi(0))+\int_0^T \left[\left(u^E,\frac{\partial\varphi}{\partial
        t}\right)- ((u^E \cdot \nabla)\, u^E, \varphi)\right]\,d \tau, 
\end{equation}
for all $\varphi \in C^\infty_0([0,T[;C^\infty(\Omega))$, with $\nabla\cdot \varphi=0$ in $\Omega$, $\varphi\cdot n=0$ on $\partial \Omega$ and all
$0\leq s<T$. 
\\
The standard argument to prove energy equality for strong solutions to~\eqref{eq:Euler} consists in choosing as test functions in
\eqref{eq:Eulerequality} the solution itself. When dealing with weak
solutions, as in our case, this procedure is not allowed, hence we
need to consider as $\varphi$ in~\eqref{eq:Eulerequality} a suitable
regularization of $u^E$. In particular we deal with a double
mollification in space and time separately, that we will denote by
$(u^E_\delta)_\varepsilon$. Here $(\cdot)_\varepsilon$ denotes the
time-mollification introduced in Subsection~\ref{moltime} and
\begin{equation}\label{uspacemol}
    u^E_\delta = \mathcal{K}^\text{div}_\delta u^E,
\end{equation}
is the divergence preserving space-mollification introduced in the
Subsection~\ref{molspace}. Note that the hypotheses of
Theorem~\eqref{thm:grad-Euler}, allows us to say that 
$u^E\in {\widetilde{V}}_{div,q}(\Omega)\subset V_{div,q}(\Omega)$,
since $\nabla u^E\in L^q(\Omega)$, and $\nabla \cdot u^E = 0$ in
$\Omega$. Moreover, being $u^{E}$ divergence free and tangential to
the boundary and by 
Lemmata~\ref{lem:smoothness}-\ref{lem:convergenceq}, $u^E$ can be
approximated by the sequence \eqref{uspacemol} that converges in
$L^q(\Omega)$ to $u^E$. A further key-point of the proof will be given
by Lemma~\ref{lem:convergeneW1q}, that gives the boundedness of the
gradient of $u_\delta^E$ in $L^q(\Omega)$ for every $\delta\in(0,1]$.

We can now give the proof of the first result of this paper.
\begin{proof}[Proof of Theorem~\ref{thm:grad-Euler}]
  Let $0<T<+\infty$ be given and 
  let $\{u_\delta^E\}$ be the sequence defined
  in~\eqref{uspacemol}. 
  For some $0<\varepsilon<T$, we choose as legitimate test function
  $\varphi=(u_\delta^E)_\varepsilon \subset
  C^\infty_0(]0,T[;C_0^\infty(\Omega))$, converging to $u^E$ in
  $L^2(0,T; L^\infty(\Omega))\cap L^p(0,T; L^q(\Omega))$. In this way
  we get the identity
\begin{equation}
  \begin{aligned}
    \label{eq:regularised_equality}
    (u^E(T), (u_\delta^E)_\varepsilon(T))&= (u_0^E, (u_\delta^E)_\varepsilon(0))
    \\
    &+\int_0^{T} \left[\left(u^E,\frac{\partial(u_\delta^E)_\varepsilon}{\partial t}\right)- ((u^E
      \cdot \nabla)\, u^E, (u_\delta^E)_\varepsilon)\right]\,d t.
  \end{aligned}
\end{equation}
Our goal is to study the previous equality and to show that in the
limit procedure it gives the energy equality for $u^E$. More precisely,
we are considering first the limit for $\delta\to 0^+$, with
$\varepsilon>0$ fixed, and then the limit as $\varepsilon\to0^+$.  It
is the passage to the limit in the nonlinear term the one which
deserves more attention. Thus, we first focus on the following term
\begin{equation*}
  \int_0^{T} ((u^E \cdot \nabla)\, u^E, (u^E_\delta)_\varepsilon)\, d t=  \int_0^{T}  \int_0^{T}
  k_\varepsilon (t-\tau) ((u^E(t) \cdot \nabla)\, u^E(t), u_\delta^E (\tau)) \,\,d \tau \,dt.
\end{equation*}
We rewrite it as the sum of three integrals as follows
\begin{equation}
  \label{eq:to-prove}
  \begin{aligned}
    \int_0^{T} ((u^E \cdot \nabla)\, u^E, (u_\delta^E)_\varepsilon) \,\,d t =& \int_0^{T} ((u^E \cdot
    \nabla)\, u^E, (u_\delta^E)_\varepsilon-(u^E)_\varepsilon) \,d t
    \\
    &+\int_0^{T} ((u^E \cdot \nabla)\, u^E, (u^E)_\varepsilon-u^E)\,
    \,d t
    \\
    &
    +\int_0^{T} ((u^E \cdot \nabla)\, u^E, u^E) \,d t
    \\
    &=: \mathbb I_1 + \mathbb I_2+\mathbb I_3. 
\end{aligned}
\end{equation}
and we will prove that: $\mathbb I_1\to 0$ as 
$\delta\to0^+$, with $\varepsilon>0$ fixed; $\mathbb I_2\to 0$ as  $\varepsilon\to0^+$; $\mathbb I_3=0$ provided
that $u$ is a weak solutions satisfying the condition of
Theorem~\ref{thm:grad-Euler}. 

\medskip

\noindent \textit{The integral $\mathbb I_3$:} As a first step we prove that 
\begin{equation} 
 \label{eq:tesi 1}
 \mathbb I_3=\int_0^{T} ((u^E \cdot \nabla)\, u^E, u^E) \,d t=0,
\end{equation}
 with $u^E$ satisfying the hypothesis of Theorem~\ref{thm:grad-Euler}.

Let $\{u_m^E\} \subset C^\infty_0(]0,T[;C^\infty(\Omega))$ be any (not
necessarily divergence-free or tangential to the boundary)
given sequence converging to $u^E$ in the space
$L^\infty(0,T;H)\cap L^p(0,T;{W}^{1,q}(\Omega))$, which exists  by density
argument. Since the field $u_m$ is smooth and $u^E$ is redefined in
such a way that $u^E(t)\in H$ for all $t\in[0,T)$, integrating by
parts, we get
\begin{equation*}
  \int_0^{T} ((u^E \cdot \nabla)\, u_m^E, u_m^E)\,d t=0.
  \end{equation*}
  Hence~\eqref{eq:tesi 1} holds true as soon as we show the following convergence
  \begin{equation*}
  \int_0^{T} ((u^E \cdot \nabla)\, u_m^E, u_m^E)\,d t \rightarrow \int_0^{T} ((u^E \cdot
  \nabla)\, u^E, u^E)\,d t.
\end{equation*}
We prove this convergence by splitting into two terms by the triangle
inequality. 
 \begin{equation}\label{eq:absolutevaluedifference}
  \begin{aligned}
    &\abs{\int_0^{T} ((u^E \cdot \nabla)\, u_m^E, u_m^E)\,d t - \int_0^{T} ((u^E \cdot
      \nabla)\, u^E, u^E)\,d t}
    \\
    &\leq \abs{\int_0^{T} ((u^E \cdot \nabla)\, u_m^E, (u_m^E-u^E))\,d t} +\abs{ \int_0^{T}
      ((u^E \cdot \nabla)\, (u_m^E- u^E), u^E)\,d t}.
  \end{aligned}
\end{equation}

Throughout the proof we will use the notation $x\lesssim y$ if there
exists a positive constant $c>0$ such that $x\leq c\,y$, where $c$
possibly depends only on $p,q,\Omega,T$, but not on the solution $u$
itself.

Let $q>\frac{9}{5}$ and
$\nabla u \in L^\frac{5q}{5q-6}(0,T,L^q(\Omega))$ and set $\widetilde{q}:=2q'$. By the
Gagliardo-Nirenberg inequality, we have
\begin{equation}
  \label{eq:Gagliardo-Nirenberg}
    \|u^E\|_{\widetilde{q}}\le C\|u^E\|_2^{\theta}\|\nabla u^E\|_q^{1-\theta},
\end{equation}
where
\begin{equation}
  \label{eq:qtilde}
    \frac{1}{\widetilde{q}}= \frac{\theta}{2}+
    (1-\theta)\bigg(\frac{1}{q}-\frac13\bigg), 
\end{equation}
and so
\begin{equation}
  \label{eq:theta}
    \theta = \frac{5q-9}{5q-6}.
\end{equation}
We estimate the second term from the right-hand side
of~\eqref{eq:absolutevaluedifference} as follows 
\begin{align*}
  \abs{ \int_0^{T} ((u^E \cdot \nabla)\,(u_m^E- u^E), u^E)\,d t} &\leq { \int_0^{T}
    \norm{u^E}_{\widetilde{q}} \norm{\nabla (u_m^E- u^E)}_q
                                                                     \norm{u^E}_{\widetilde{q}}\,d t} 
  \\
  \lesssim &{ \int_0^{T} \|u^E\|_2^{2\theta}\|\nabla
             u^E\|_q^{2(1-\theta)}\norm{\nabla (u_m^E- u^E)}_q \,d t} 
  \\
  \lesssim &\|u^E\|_{\infty,2}^{2\theta}\|\nabla
             u^E\|_{p,q}^{2(1-\theta)}\norm{\nabla (u_m^E-
             u^E)}_{p,q}, 
\end{align*}
where in the first line we have applied the H\"older inequality with
the three conjugate exponents $\widetilde{q},q,\widetilde{q}$, in the
second line the Gagliardo-Nirenberg inequality and in the third line
the H\"older inequality in the time variable with conjugate exponents
$r,s$ such that $2(1-\theta)r=p$ and $s=p$, and $2\theta y=p$, which
gives $p= \frac{5q}{5q-6}$. Since $\nabla u^E_m\to \nabla u^E$ in
$L^p(L^q(\Omega))$, we have that
$$\int_0^{T} ((u^E \cdot \nabla)\,(u_m^E-u^E), u^E)\,d t\rightarrow
0,\qquad \text{as}\,\,\,\, m\to +\infty.$$ Let us stress that when
$q=9/5$, i.e. $\nabla u^E \in L^3(L^\frac{9}{5})$, the computation is
similar. In this case, equations~\eqref{eq:qtilde}-\eqref{eq:theta}
show that $\theta=0$ and $\widetilde{q}=q^*$ is the Sobolev exponent,
so that the Gagliardo-Nirenberg
inequality~\eqref{eq:Gagliardo-Nirenberg} becomes the classical the
Sobolev embedding with $q=9/5$. Therefore, using the H\"older
inequality in space with the three exponents $(9/5)^*$, $9/5$,
$(9/5)^*$, the Sobolev embedding and the H\"older inequality in time
with conjugate exponents $3$ and $3/2$, we get
\begin{align*}
  \abs{ \int_0^{T} ((u^E \cdot \nabla)\,(u_m^E- u^E), u^E)\,d t}
  &\lesssim\|\nabla u^E\|_{3,\frac{9}{5}}^2\norm{\nabla (u_m^E- u^E)}_{3,\frac{9}{5}},
\end{align*}
arriving at the same conclusion as in the previous case.

Analogously, we estimate the first term of the right-hand side
of~\eqref{eq:absolutevaluedifference}, when $q>9/5$ 
\begin{align*}
  &\abs{ \int_0^{T} ((u^E \cdot \nabla)\,u_m^E ,(u_m^E-u^E))\,d t}
  \\
  &\leq \abs{ \int_0^{T} \norm{u^E}_{\widetilde{q}} \norm{\nabla u_m^E}_q
    \norm{u_m^E-u^E}_{\widetilde{q}}\,d t}
  \\
  &\lesssim\abs{ \int_0^{T} \norm{u^E}_{2}^{\theta} \norm{
      \nabla u^E}_{q}^{(1-\theta)}\norm{\nabla u_m^E}_q \norm{u_m^E-u^E}_{2}^{\theta} \norm{
      \nabla(u_m^E-u^E)}_{q}^{(1-\theta)}\,d t}
      \\
  &\lesssim \norm{u^E}_{\infty,2}^{\theta}\norm{u_m^E-u^E}_{\infty,2}^{\theta}\abs{ \int_0^{T}  \norm{
      \nabla u^E}_{q}^{(1-\theta)}\norm{\nabla u_m^E}_q  \norm{
      \nabla(u_m^E-u^E)}_{q}^{(1-\theta)}\,d t}
  \\
  &\lesssim \norm{u^E}_{\infty,2}^{\theta}\norm{u_m^E-u^E}_{\infty,2}^{\theta}  \norm{
      \nabla u^E}_{p,q}^{(1-\theta)}\norm{\nabla u_m^E}_{p,q}  \norm{
      \nabla(u_m^E-u^E)}_{p,q}^{(1-\theta)} 
\end{align*}
where again we have used the Gagliardo-Nirenberg inequality in the
third line and the H\"older inequality in the time variable with
conjugate exponents $x(1-\theta)=p$, $y=p$ and $z(1-\theta)=p$. For
$m\to +\infty$ we have the desired convergence.\\

When $q=9/5$, we argue as before.
By using the H\"older
inequality in space with the three exponents $(9/5)^*$, $9/5$,
$(9/5)^*$, the Sobolev embedding and H\"older inequality in time with
the three conjugate exponents $3$, $3$, $3$, we get
\begin{align*}
  &\abs{ \int_0^{T} ((u^E \cdot \nabla)\,u_m^E ,(u_m^E-u^E))\,d t}
  &\lesssim \norm{
      \nabla u^E}_{3,\frac{9}{5}}\norm{\nabla u_m^E}_{3,\frac95}  \norm{
      \nabla(u_m^E-u^E)}_{3,\frac95}.
\end{align*}
Even in this case we have the desired convergence when
$m \to +\infty$. Collecting the above estimates we conclude the proof
of~\eqref{eq:tesi 1}.

\medskip

\noindent\textit{The integral $\mathbb I_1$:} Next we
show that the first integral from the right hand
side of~\eqref{eq:to-prove} converges to zero.

 %
When $q>9/5$, applying again the H\"older, the Gagliardo-Nirenberg and
the Sobolev inequalities we get
\begin{align*}
  &\abs{\mathbb I_1}=\abs{ \int_0^{T} ((u^E \cdot \nabla)\,u^E,
    (u_\delta^E)_\varepsilon-(u^E)_\varepsilon)\,d t}  
  \\
  &\leq {\int_0^{T}\norm{u^E}_{\widetilde{q}} \norm{\nabla u^E}_q
    \norm{(u_\delta^E)_\varepsilon-(u^E)_\varepsilon}_{\widetilde{q}}\,d t} 
  \\
  &\lesssim { \int_0^{T} \norm{u^E}_{2}^{\theta} \norm{\nabla u^E}_{q}^{1-\theta} 
    \norm{\nabla u^E}_q  \norm{(u_\delta^E)_\varepsilon-(u^E)_\varepsilon}_{2}^{\theta}
    \norm{\nabla((u_\delta^E)_\varepsilon-(u^E)_\varepsilon)}_{q}^{1-\theta} \,d t}
  \\
  &\lesssim \norm{u^E}_{\infty,2}^{\theta}\norm{(u_\delta^E)_\varepsilon-(u^E)_\varepsilon}_{\infty,2}^{\theta}{ \int_0^{T}  \norm{\nabla u^E}_{q}^{2-\theta} \norm{\nabla((u_\delta^E)_\varepsilon-(u^E)_\varepsilon)}^{1-\theta}_{q} \,d t}
  \\
  &\lesssim \norm{u^E}_{\infty,2}^{\theta} \norm{(u_\delta^E)_\varepsilon-(u^E)_\varepsilon}_{\infty,2}^{\theta} \norm{\nabla u^E}_{p,q}\norm{\nabla((u_\delta^E)_\varepsilon-(u^E)_\varepsilon)}_{p,q}
\end{align*}
where in the last line we have applied the H\"older inequality in the
time variable with conjugate exponents $x,y$, such that
$x(2-\theta)=p$ and $y(1-\theta)=p$, which gives
$p=\frac{5q}{5q-6}$. Now, Corollary~\ref{cor:LinftyLq} ensures that
\begin{equation*}
  \lim_{\delta\rightarrow 0^+} \norm{(u_\delta^E)_\varepsilon-(u^E)_\varepsilon}_{\infty,2}=0, 
\end{equation*}
and, moreover by Corollary~\ref{cor:boundednessgradtimespace} we also
have  
  $\norm{\nabla((u_\delta^E)_\varepsilon-(u^E)_\varepsilon)}_{p,q}\leq C $ 
which, together with the previous estimate, implies that
\begin{equation*}
  \lim_{\delta \to 0^+} \mathbb I_1=\lim_{\delta\to0^+}\int_0^{T}
  ((u^E \cdot \nabla)\,u^E,
  (u^E_\delta)_\varepsilon-(u^E)_\varepsilon)\,d t =0, 
\end{equation*}
holds for each fixed $\varepsilon>0$. The case $q=9/5$ can be treated
in an easier way than the previous one, since by a simple application
of H\"older and Sobolev inequalities in space and the H\"older
inequality with conjugate exponents $3/2$ and $3$ in time, we can
write
\begin{align*}
  \abs{\mathbb I_1}
  &\lesssim \norm{\nabla u^E}_{3,\frac{9}{5}}^2 
     \norm{(u_\delta^E)_\varepsilon-(u^E)_\varepsilon}_{3,(\frac{9}{5})^*},
\end{align*}
obtaining the desired convergence.

Finally, the proof of the convergence of $\mathbb I_2\to0$  as
$\varepsilon\to 0^+$ follows by  similar estimates.

Eventually we conclude the passage to the limit in~\eqref{eq:regularised_equality} as $\delta$ and $\varepsilon$ go to
zero. All the previous convergence results allow to conclude that the
convective term tends to zero, that is
\begin{equation*}
\lim_{\varepsilon\to0^+}\lim_{\delta\to 0^+}\int_0^{T} ((u^E \cdot \nabla)\,u^E,
(u^E_\delta)_\varepsilon)\,d t=0. 
\end{equation*}
The term
involving the time derivative of $k_\varepsilon$ vanishes identically,
i.e.
$$
\int_0^{T}
\left(u^E,\frac{\partial(u^E_\delta)_\varepsilon}{\partial
    t}\right)=0,
$$ since
$k_\varepsilon$ is even. The mollification in space and time allow to directly prove 
\begin{equation*}
       \begin{split}
       &\lim_{\varepsilon\to 0^+}\lim_{\delta\to
         0^+}(u^E_0,(u^E_\delta)_\varepsilon(0))=\frac{\|u^E_0\|^2}{2},\\
         &\lim_{\varepsilon\to 0^+}\lim_{\delta\to
         0^+}(u^E(T),(u^E_\delta)_\varepsilon(T))=\frac{\|u^E(T)\|^2}{2},     \end{split} 
\end{equation*}
under the assumption of weak-$L^2$-continuity, see
Galdi~\cite{Gal2000a}.
In conclusion, from \eqref{eq:regularised_equality}, we get

\begin{equation*}
  \norm{u^E(T)}^2= \norm{u^E_0}^2,
\end{equation*}
which is~\eqref{eq:energy-equality} for $t=T$. We notice that the energy equality is still valid for almost every $t_0\in (0,T)$, $0<\varepsilon<t_0$, just by considering a redefinition of the mollification in time as an integral between $0$ and $t_0$. The proof
of Theorem~\ref{thm:grad-Euler} is now complete.
\end{proof}
\subsection{Energy equality for Euler Beltrami fields}
After having finished the proof of the criterion for energy
conservation about the gradient, we can pass to prove to a criterion
for energy conservation, using the relation between the vorticity and the
velocity for solutions with the special geometrical
constraint~\eqref{eq:Beltrami}. This should be also compared with the
results in~\cite{Der2020}, where an ``analytic'' combination of the
two quantities was considered, and with the results
from~\cite{BS2024}, for fractional spaces (both references deal with
the space periodic setting).
\begin{proof}[Proof of Theorem~\ref{thm:BeltramifieldthmEuler}]
  Let $u^E\in L^\infty(0,T;L^2(\Omega))$ be a weak solution to
  problem~\eqref{eq:Euler} which is a Beltrami
  field~\eqref{eq:Beltrami} and let
  $\lambda \in L^\alpha(0,T;L^\beta(\Omega))$.  Direct applications of
  the H\"older inequality in space and time show that
  \begin{equation}
    \label{eq:iteration-omega-1}
    \omega^E \in L^\alpha(0,T;L^\frac{2\beta}{2+\beta}(\Omega)),
\end{equation}
which in turn implies that $\nabla u^E$ is in the same Bochner space
(again by using the results from~\cite{vWah1992}). It is now
convenient to make a change of notation and let us determine $q$ is
terms of $(\alpha,\beta)$ as follows
\begin{equation*}
        \alpha = \frac{5q}{5q-6}\qquad\text{and}\qquad
        \frac{2\beta}{2+\beta}= q.
\end{equation*}
in such a way that~\eqref{eq:iteration-omega-1} is equivalent to the
condition~\eqref{eq:grad-Euler}  on $\nabla u^E$.
(This is possible if $(\alpha,\beta)$ are as in
Theorem~\ref{thm:BeltramifieldthmEuler}) 
Equivalently, the statement of the theorem can be rewritten, in terms of $q$, as
follows
$$
\lambda \in
L^\frac{5q}{5q-6}(0,T,L^\frac{2q}{2-q}(\Omega))\qquad\text{ with
}\quad\frac{6}{5}< q<2, 
$$
and different choices of the leading parameter can be used to a better
interpretation of the results.
Hence, if $q\in[9/5,2[$ (or equivalently $\alpha\in]5/2,3]$) by
Theorem~\ref{thm:grad-Euler} we have conservation of energy.

Let now $q\in (6/5,9/5)$ and let us show that if we iterate further
the procedure, we still have conservation of energy after a number of
suitable iterations. Indeed, notice that since
$\nabla u^E \in L^\frac{5q}{5q-6}(0,T;L^q(\Omega))$, then
$u^E \in L^\frac{5q}{5q-6}(0,T;W^{1,q}(\Omega))$, and this follows by
the Poincar\`e inequality which is valid even for tangential vector
fields, in the case of a bounded domain, see
Galdi~\cite{Gal2011}. Then, the Sobolev embedding theorem implies also
that $u^E \in L^\frac{5q}{5q-6}(0,T;L^{q^*}(\Omega)) $, where
$q^{*}= \frac{3q}{3-q}$. Since $u^E$ is a Beltrami field, then we have
again (by further applications of H\"older inequality) that
\begin{equation*}
    \nabla u^E \in L^\frac{5q}{2(5q-6)}(0,T;L^\frac{6q}{12-5q}(\Omega)).
  \end{equation*}
  If $q_1:=\frac{6q}{12-5q}\geq 9/5$ (which corresponds to
  $q\in[36/25,9/5[$) we end the proof due to the fact that
  $p_1=\frac{5q}{2(5q-6)}=\frac{5 q_1}{5q_1-6}$ and this fits again
  with the hypotheses of Theorem~\eqref{eq:grad-Euler}.

  On the other hand, if $q_1<\frac{9}{5}$ this is not enough to end
  the proof but we can further iterate the same argument. This leads to
  defining the following two sequences of indices
\begin{equation*}
  p_n := \frac{5q}{(n+1)(5q-6)}\qquad \text{and}\qquad
  q_n := \frac{6q}{6(n+1)-5n q},\qquad \text{for }n\ge 1, 
  \end{equation*}
 and  we observe that the iterative procedure (and the step $k\mapsto k+1$
  is well defined if $p_{k}\geq1+1/(k+1)$ and $q_{k}<3$, for $k\in\N$)
  implies that
\begin{equation*}
  \nabla u^E \in L^{p_n}(0,T;L^{q_n}(\Omega))\qquad\text{with}\quad 
  p_n=
  \frac{5q_n}{5q_n-6}. 
\end{equation*}
Hence, in view of Theorem~\ref{thm:grad-Euler}, the condition $\nabla
u^E \in L^{p_n}(0,T;L^{q_n}(\Omega))$ implies conservation of energy,
as soon as we can find $n_{0}\in\N$ such that that one term of
$q_{n_{0}}$ of the strictly increasing sequence
$\{q_{n}\}$ is such that $q_{n_{0}}\ge 9/5$. This condition is equivalent (in
terms of $q\in]6/5,36/25[$) to the following requirement
\begin{equation*}
\frac{6}{5}<\frac{18(n+1)}{5(3n +2)}\le
q<\frac{6}{5}+\frac{6}{5n}\leq\frac{9}{5}.   
\end{equation*}
%
The lower bound for $q$ is strictly decreasing to $6/5$ as $n\to +\infty$.
In this way, for any $q\in (6/5,36/25)$, it is sufficient to choose
$n_0\in \mathbb N$ such that $n_0\geq\frac{18-10q}{15q-18}$ to ensure
that $q_{n_0}\ge 9/5$. Thus, after $n_0\in\N$ steps of the iteration,
Theorem~\ref{thm:grad-Euler} can be applied. Collecting all the
results this allows us to conclude the proof in the range
$q \in (6/5,2)$.

The case $q=2$, which corresponds to $\alpha=5/2$, can be studied
separately. Indeed if
$\lambda \in L^\frac{5}{2}(0,T;L^\infty(\Omega))$, then we can stop at
the first iteration, since we have %
%
$  \nabla u^E \in L^{\frac{5}{2}}(0,T;L^2(\Omega))$, %
%
which is~\eqref{eq:grad-Euler} for $q=2$.
\end{proof}
\begin{rem}
  The iterative boot-strap argument we used does not ``improve'' the
  known regularity of the solution $u^{E}$, in terms of scaling. This
  means that at each step we get that
  $\nabla u^{E}\in L^{5q/(5q-6)}(0,T;L^{q}(\Omega))$, with different
  values of $q$; the iteration is needed just to ensure that
  $q\geq9/5$, to enter the range of validity of
  Theorem~\ref{thm:grad-Euler}. The situation will be rather different
  in the case of the Navier-Stokes equations, as we will see in the
  next section.
\end{rem}

\section{On the energy equality and regularity for
  NSE weak (Beltrami) solutions}\label{sec:BeltramiNSE}
In this section we study the energy conservation and regularity for a
Leray-Hopf weak solution to the Navier-Stokes
equations~\eqref{eq:NSE}, which is also a Beltrami field, proving
Theorems~\ref{thm:Beltrami-NSE}-\ref{thm:regularityNSE}

\begin{proof}[Proof of Theorem~\ref{thm:Beltrami-NSE}]
  The proof is very similar to that of
  Theorem~\ref{thm:BeltramifieldthmEuler} but one relevant difference
  --this also justifies the slight change of notation-- is that we
  have to check also that the iteration process does not imply that
  the solution enters (by a bootstrap argument) into a class of
  regularity, for which the energy equality holds trivially. This is
  particularly delicate for the NSE: contrary to the Euler equations
  where practically only the Beale-Kato-Majda criterion holds, various
  less restrictive scaling invariant conditions are known in the
  viscous case, as for instance~\eqref{eq:scaling-nabla} that we will
  employ.

  Hence, in this case the task will be not only to show that the
  gradient belongs to a space with the exponents as in
  Theorem~\ref{thm:Berselli-Chiodaroli}, but also to show that they do
  not reach the class as in~\eqref{eq:scaling-nabla}. Note that the
  iteration we handled in the previous theorem was needed just to
  reach an exponent $q$ in the space variables larger or equal than
  $9/5$. The exponents $p_{n},q_{n}$ defined in the proof of the
  previous theorem satisfy conditions as in
  Theorem~\ref{thm:Berselli-Chiodaroli}, which is out of the
  regularity class, for all $n\leq n_{0}$. Hence if we perform
  iterations up to first $n\in\N$ such that $q_{n}>3$, the gradient
  still belongs to~\eqref{eq:grad-Euler} and hence not
  to~\eqref{eq:scaling-nabla}. 

  In this case we have also to see what happens further iterating the
  sequence: a further iteration shows that the velocity belongs to
  $L^{\infty}(\Omega)$, possibly drastically changing the known
  scaling of the solution.

 Let then  $u$ be a Leray-Hopf weak solution such that 
  $  \omega(x,t)= \lambda(x,t)u(x,t)$, %
with $\lambda \in L^\alpha(0,T;L^\beta(\Omega))$, for some
$\alpha,\beta \ge 1$ to be determined in order to obtain conservation
of energy.
Since $u\in L^\infty(0,T;L^2(\Omega))$, by
using~\eqref{eq:iteration-omega-1} we get 
$ \nabla u \in L^\alpha(0,T;L^\frac{2\beta}{2+\beta}(\Omega))$. Before
introducing the iteration procedure, we rename the time and space
exponents as follows
\begin{equation*}
    p_1 = \alpha, \qquad \text{and}\qquad q_1 = \frac{2\beta}{2+\beta}.
\end{equation*}
In order to apply Theorem~\ref{thm:Berselli-Chiodaroli}, we need to
require that 
\begin{equation*}
    \begin{cases}
      q_1>\frac{9}{5},
      \\
        p_1 = \frac{5q_1}{5q_1-6},
    \end{cases}\iff
    \begin{cases}
      \beta>18,
      \\
        \beta = \frac{6\alpha}{2\alpha-5},
    \end{cases}
\end{equation*}
which is equivalent to require that 
$        \beta = \frac{6\alpha}{2\alpha-5}$ for
      $\frac{5}{2}<\alpha <3$,
and so 
\begin{equation}
  \label{p1q1}
    p_1 = \alpha, \qquad q_1 = \frac{6\alpha}{5\alpha-5}
    \qquad\text{with}\quad 
    \frac{5}{2}<\alpha <3. 
\end{equation}
These requirements guarantee the conservation of energy. Moreover, we
stress again the fact that the solution does not belong to the regularity
class~\eqref{eq:scaling-nabla}; indeed from~\eqref{p1q1} we have (and
recall that 
at this step $q_1<3$)
\begin{equation*}
  \frac{2}{p_1}+\frac{3}{q_1}=\frac{5}{2}-\frac{1}{2\alpha}=
  2+\bigg(\frac{1}{2}-\frac{1}{2\alpha}\bigg)> 2.
\end{equation*}
If instead $ q_1 = \frac{2\beta}{2+\beta} < \frac{9}{5}$, which  means
$    \beta <18$ (or equivalently $\alpha > 3$), 
then, as in the proof of the previous
Theorem~\ref{thm:BeltramifieldthmEuler}, we need to iterate the
procedure. From $\nabla u \in L^{p_1}(0,T;L^{q_1}(\Omega))$, we get
$u \in L^{p_1}(0,T;L^{q_1^*}(\Omega))$ by the Sobolev embedding
theorem, with
  $  q_1^*= \big(\frac{2\beta}{2+\beta}\big)^* =
  \frac{6\beta}{6+\beta}$. 
Since $u$ is a Beltrami field, then it follows that 
\begin{equation*}
    \nabla u \in
    L^{\frac{\alpha}{2}}(0,T;L^\frac{6\beta}{12+\beta}(\Omega)). 
\end{equation*}
Here, we call
  $  p_2 = \frac{\alpha}{2}$ and $q_2 = \frac{6\beta}{12+\beta}$ and 
to have conservation of energy we need to require that
\begin{equation*}
    \begin{cases}
      q_2 > \frac95,
      \\
        p_2= \frac{5q_2}{5q_2-6},
        \end{cases} \quad\iff\quad 
        \begin{cases}
          \beta > \frac{36}{7},
          \\
            \beta = \frac{6\alpha}{2\alpha-5},
        \end{cases}
\end{equation*}
that can be satisfied if $3<\alpha < 6$. In particular, we can rewrite
everything in terms of $\alpha$ as $\nabla u\in
L^{\frac{\alpha}{2}}(0,T;L^{\frac{6\alpha}{5\alpha-10}}(\Omega))$, 
and even in this case we have conservation of energy, but we do not fall
in the regularity class~\eqref{eq:scaling-nabla} since
\begin{equation*}
  \frac{2}{p_2}+\frac{3}{q_2}=2 +
  \bigg(\frac{1}{2}-\frac{1}{\alpha}\bigg)>2. 
\end{equation*}
If $q_2 < \frac{9}{5}$, then we need to iterate again. So proceeding
as before we can construct two sequences of exponents (written in
terms of $\alpha$ as follows)
\begin{equation*}
    p_n = \frac{\alpha}{n}, \qquad q_n = \frac{6\alpha}{5\alpha-5n},
\end{equation*}
with $\alpha > n$ and such that at the $n$-th step we have
\begin{equation}
  \label{nablapnqn}
    \nabla u \in L^{p_n}(0,T;L^{q_n}(\Omega)), \qquad n\ge 3,
\end{equation}
where note that 
  $  p_n = \frac{5q_n}{5q_n-6}$,
  cf. Thm.~\ref{thm:Berselli-Chiodaroli}. 
In particular if at the $(n-1)$-th step we have $q_{n-1}<\frac{9}{5}$,
then to have conservation of energy we need to require that
$q_n>\frac{9}{5}$, that gives
%
  $  \alpha < 3n$. %
%
Even in this case, to have a non-trivial result we need to check if we
do not fall in the regularity class. Indeed since $\alpha > n$ we have
\begin{equation*}
    \frac{2}{p_n}+\frac{3}{q_n}=
    2+\bigg(\frac12-\frac{n}{2\alpha}\bigg) > 2. 
\end{equation*}
The previous iteration argument shows that for any fixed
$\alpha > \frac{5}{2}$ and $\beta = \frac{6\alpha}{2\alpha-5}$, in
order to have conservation of energy we need to iterate the algorithm
a number $n_0\in\N$ of times, such that
\begin{equation*}
    \frac{\alpha}{3}< n_0 < [\alpha],
\end{equation*}
where $[\cdot]$ denotes the integer-part function.

To conclude the proof we need to be careful about $q_n$ being 
smaller or larger than $3$. Indeed, as long as $q_n<3$ we keep on
iterating and constructing the sequences of exponents as we showed
previously, which are the same also of
Theorem~\ref{thm:BeltramifieldthmEuler}, but written with a more
convenient notation.  On the other hand, if $q_n>3$, even if we are
already in a class of energy conservation we can still further iterate
the process to see if there is a further improvement in the
regularity (beyond energy conservation). In the case $q_n>3$  Sobolev
embedding theorem applied to~\eqref{nablapnqn} implies that
   $ u \in L^{p_n}(0,T;L^\infty(\Omega))$, %
  %
which consequently proves that
\begin{equation*}
    \nabla u \in L^{\frac{\alpha}{n+1}}(0,T;L^\beta(\Omega)).
  \end{equation*}
  \begin{rem}
    Note that the space integrability for the gradient never goes
    beyond $\beta$, being obtained through multiplication by $\lambda$
    and the H\"older inequality.
  \end{rem}
In this last case we get (in terms of $\alpha$)
\begin{equation}
  \label{pn1qn1}
    p_{n+1}= \frac{\alpha}{n+1},\qquad\text{and}\qquad q_{n+1} =
    \beta=\frac{6\alpha}{2\alpha-5}. 
\end{equation}
Even though we have conservation of energy (since solution was already
in a class as in Theorem~\ref{thm:Berselli-Chiodaroli}), in this case
we still need to check if~\eqref{eq:scaling-nabla} holds. The
condition $q_n>3$ is equivalent to have
\begin{equation}
  \label{nge35alpha}
    n> \frac{3}{5}\alpha,
\end{equation}
and so, considering~\eqref{pn1qn1},  we get
\begin{equation*}
    \frac{2}{p_{n+1}}+\frac{3}{q_{n+1}} =\frac{4(n+1)-5}{2\alpha}+1,
\end{equation*}
and the right-hand side is strictly greater than $2$ if and only if
\begin{equation*}
    \frac{4(n+1)-5}{2\alpha}>1 \iff n > \frac{\alpha}{2}+\frac{1}{4}.
\end{equation*}
But, since~\eqref{nge35alpha} is satisfied, we ask whether  %
$   \frac{3}{5}\alpha\ge \frac{\alpha}{2}+\frac{1}{4}$,
and this gives $\alpha \ge \frac{5}{2}$. This shows that even in the
case of $q_n>3$ we do not fall in a regularity class. Note that a
further iteration will not improve the space exponent, but will reduce
the one in the time variable, hence not improving the known regularity
of $u$.  Hence, the proof is concluded in the case of
$\alpha > \frac{5}{2}$.

 The case $\alpha = \frac{5}{2}$ follows the same argument as in
 the proof of Theorem~\ref{thm:BeltramifieldthmEuler}. Indeed if
 $\lambda \in L^\frac{5}{2}(0,T;L^\infty(\Omega))$, then we can stop
 at the first iteration, since we have
 \begin{equation*}
         \nabla u \in L^{\frac{5}{2}}(0,T;L^2(\Omega)),
     \end{equation*}
     which is not a regularity class but  falls immediately in the
     case $(ii)$ of Theorem~\ref{thm:Berselli-Chiodaroli}. Further
     iterations will not improve the regularity, cf. also with 
     Proposition~\ref{prop:Elementary}. 
 \end{proof}
 Since we identified classes of solutions for which the energy
 conservation holds true but these weak solutions are not strong
 (hence regular), we now prove Theorem~\ref{thm:regularityNSE} which
 describes natural assumptions on $\lambda$ implying regularity of
 solutions (hence a fortiori also uniqueness and energy conservation).
 \begin{proof}[Proof of Theorem~\ref{thm:regularityNSE}]
   Let $u$ be a Leray-Hopf solution to problem~\eqref{eq:NSE},
   which is also a Beltrami field. Let
   $\lambda \in L^\alpha(0,T;L^\beta(\Omega))$. As before we get as first step
   $ \nabla u \in L^\alpha (0,T;L^\frac{2\beta}{2+\beta})$ and this
   implies that $u$ is a strong solution if
\begin{equation*}
  \frac{2}{\alpha}+\frac{3(2+\beta)}{2\beta} = 2\quad \iff\quad  \frac{2}{\alpha}+\frac{3}{\beta}= \frac{1}{2}, 
\end{equation*}
with $ \frac{3}{2}< \frac{2\beta}{2+\beta}$, which gives $\beta >
6$ and within this range
$ \alpha = \frac{4\beta}{\beta-6}$.
Let us notice that $ \frac{2\beta}{2+\beta}<3$, and this tells us that
we can iterate the procedure. As in the step $2$ of previous theorem
we get
   $ \nabla u \in
    L^{\frac{\alpha}{2}}(0,T;L^\frac{6\beta}{12+\beta}(\Omega))$. 
and the scaling invariant regularity condition connected to the above integrability exponents reads as
\begin{equation}
  \label{eq:regularitysecondstep}
    \frac{2}{\alpha}+\frac{3}{\beta}= \frac{3}{4}.
\end{equation}
Now we have to distinguish two cases.\\
\textit{i)} If $\beta > 12$, then the space exponent is
$q_2=\frac{6\beta}{12+\beta}>3$. In this case the Sobolev embedding
yields $ u \in L^{\frac{\alpha}{2}}(0,T;L^\infty(\Omega))$, so that
the next step of the iteration gives
  $  \nabla u \in L^{\frac{\alpha}{3}}(0,T;L^\beta(\Omega))$,
which corresponds to the following regularity condition
\begin{equation}\label{eq:regularitylinfty}
    \frac{2}{\alpha}+\frac{3}{\beta}= 2-\frac{4}{\alpha}.
\end{equation}
In order to optimize the choice of $\alpha$ and $\beta$, we need to compare the two conditions~\eqref{eq:regularitysecondstep} and~\eqref{eq:regularitylinfty}. To this aim we express $\alpha$ as a function of $\beta$, getting respectively
\begin{equation*}
    \alpha=\frac{8\beta}{3\beta-12}\qquad \text{and}\qquad \widetilde{\alpha} = \frac{6\beta}{2\beta-3}.
\end{equation*}
It is easy to check that $\alpha<\widetilde{\alpha}$ if and only if
$\beta > 24$. Thus if $\beta > 24$, we choose $\alpha$ such that the
condition~\eqref{eq:regularitysecondstep} holds, that corresponds
to~\eqref{eq:regularityLnRn}$_2$ when $\overline{n}=1$. If
$\beta \in (12,24]$, $\widetilde{\alpha}$ is the optimal exponent and
the condition~\eqref{eq:regularityLnRn}$_1$ (in the case of
$\overline{n}=1$) is obtained by~\eqref{eq:regularitylinfty} with the
choice of $\widetilde{\alpha}=\widetilde{\alpha}(12) = 24/7$. Notice
that condition~\eqref{eq:regularityLnRn}$_1$, for $\overline{n}=1$, implies
condition~\eqref{eq:regularitylinfty} for any choice of
$\widetilde{\alpha}=\widetilde{\alpha}(\beta)$, with $\beta \in
(12,24]$.\\
\textit{ii)} Let now $\beta\in (3,12)$. In this case $q_2<3$ and we
can keep iterating. As we have already done in the previous theorem,
we find two sequences of exponents
\begin{equation*}
  p_n = \frac{\alpha}{n},\qquad q_n = \frac{6\beta}{6n-(2n-5)\beta}
  \qquad n\ge 2, 
\end{equation*}
such that $ \nabla u \in L^{p_n}(0,T;L^{q_n}(\Omega))$.  The solution
is regular if
%
\begin{equation*}
    \frac{2}{p_n}+\frac{3}{q_n}= 2\  \iff\
    \frac{2}{\alpha}+\frac{3}{\beta}=1-\frac{1}{2n},
\end{equation*}
which gives
\begin{equation}
  \label{eq:alphan}
    \alpha_n = \frac{4n\beta}{2n(\beta-3)-\beta}.
\end{equation}
In order to go further in the iteration process we have to impose
$ \frac{3}{2}< q_n < 3$, that reads as
\begin{equation*}
    \frac{6n}{2n-1}<\beta < \frac{6n}{2n-3}.
\end{equation*}
Let us now suppose that at the $(n+1)-$step, we have $q_{n+1}>3$, that
is more precisely 
\begin{equation}
  \label{eq:qn3qn+1}
      q_n<3
\qquad\text{and}\qquad 
        q_{n+1}>3,
 \end{equation}
which is equivalent to 
%
%
\begin{equation*}
  \beta\in  I_n := \left( \frac{6(n+1)}{2n-1},\frac{6n}{2n-3}\right).
\end{equation*}
In particular, restating the previous inequality in terms of $n$, we
fix the number of iterations we can make. Indeed we have
\begin{equation*}
    \frac{\beta+6}{2\beta-6}< n < \frac{3\beta}{2\beta-6},
\end{equation*}
and since the size the gap is 
$    \frac{3\beta}{2\beta-6}-\frac{\beta+6}{2\beta-6}=1$, %
then, for any fixed $\beta \in (3,12)$, there exists a unique
$n\ge 2$, such that~\eqref{eq:qn3qn+1} occurs.

Now, since $q_{n+1}>3$, the condition
$\nabla u \in L^{p_{n+1}}(0,T;L^{q_{n+1}}(\Omega))$ implies by the
Morrey-Sobolev embedding that  
$u \in L^{p_{n+1}}(0,T;L^\infty(\Omega))$ and, iterating once more,
\begin{equation*}
    \nabla u \in L^{p_{n+2}}(0,T;L^\beta(\Omega)).
\end{equation*}
In this case we have regularity of the solution if 
\begin{equation}
  \label{eq:regularityn+2}
    \frac{2(n+2)}{\alpha}+\frac{3}{\beta}=2 \iff \frac{2}{\alpha}+\frac{3}{\beta}= 2-\frac{2(n+1)}{\alpha}.
  \end{equation}
  In particular, we can express the integrability exponent of
  $\lambda$ with respect to time in terms of $\beta$ as follows:
\begin{equation*}
    \widetilde{\alpha}_{n+2}= \frac{2(n+2)\beta}{2\beta-3}.
\end{equation*}
We stress the fact that at this point
we can stop iterating, since $\nabla u$ would end up in a more
regular class.

To conclude, we need to understand what is the best choice between
$\alpha_{n+1}$, given in~\eqref{eq:alphan}, and
$\widetilde{\alpha}_{n+2}$ and, to start, we observe that
\begin{equation*}
    \alpha_{n+1}> \tilde{\alpha}_{n+2}\qquad\iff\qquad
    \beta < \frac{6(n+1)^2}{2n^2+n-2} \in I_n.
\end{equation*}
We  already defined $L_{n}$ and $R_n$ in~\eqref{eq:LnRn} %
and the above argument tells us that we need to choose:
\begin{equation*}
    \begin{cases}
      \widetilde{\alpha}_{n+2}&\text{when}\,\,\,\beta \in L_n
      \\
        \alpha_{n+1} &\text{when}\,\,\,\beta \in R_n.
    \end{cases}
\end{equation*}

\begin{figure}[H]
\centering
\includegraphics[scale=0.1
]{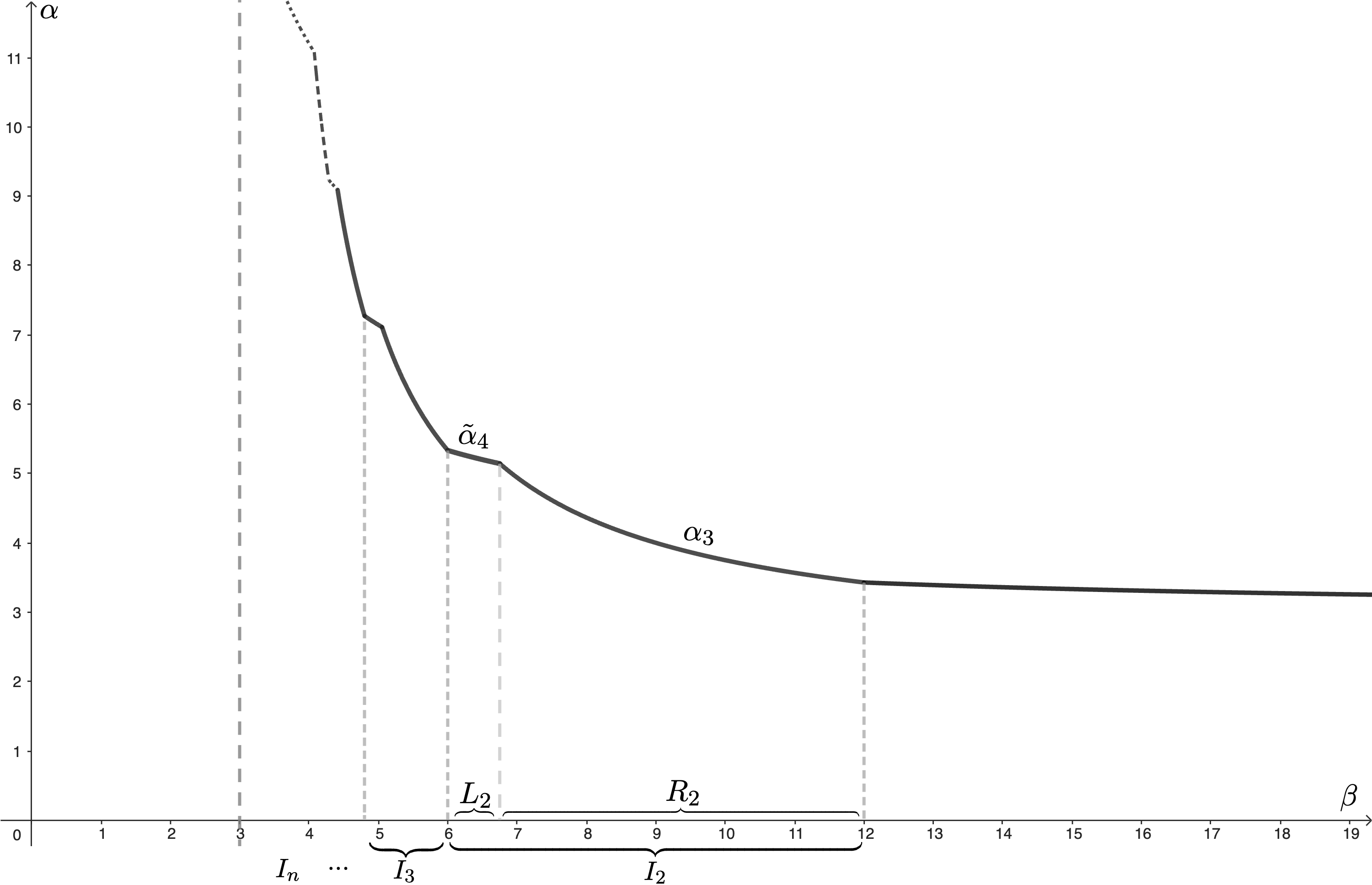}
\caption{Graphs of $\widetilde{\alpha}_{n+2}$ and $\alpha_{n+1}$ as functions of $\beta$ in the intervals $L_n$ and $R_n$.}
\end{figure}
The choice of $\alpha_{n+1}$ or $\widetilde{\alpha}_{n+2}$ corresponds
respectively to the choice of $\alpha$ satisfying
\begin{equation}
  \label{eq:regularitiesalpha}
  \frac{2}{\alpha}+\frac{3}{\beta}= 1-\frac{1}{2(n+1)}
  \qquad
    \text{or}\qquad\frac{2}{\alpha}+\frac{3}{\beta}=
    2-\frac{2(n+1)}{\alpha}, 
\end{equation}
when $\beta$ belongs to $L_{n}$ or $R_{n}$, respectively (note that $\alpha>n+1$).

Let us notice that if we evaluate $\widetilde{\alpha}_{n+2}$ in $\beta_0 = \frac{6(n+1)}{2n-1}$ (which is the left endpoint of $L_n$), we get
\begin{equation*}
    \widetilde{\alpha}_{n+2}(\beta_0)= \frac{4(n+1)(n+2)}{2n+5},
\end{equation*}
and substituting the latter in the left-hand side of the second
equation in~\eqref{eq:regularitiesalpha} we have
\begin{equation}
  \label{eq:2alpha3betan2}
  \frac{2}{\alpha}+\frac{3}{\beta}= 1-\frac{1}{2(n+2)}, 
\end{equation}
showing that $\alpha$ is continuous also between the intervals $I_n$
and $I_{n+1}$. In particular, condition~\eqref{eq:2alpha3betan2}
implies the one contained in~\eqref{eq:regularityn+2}.

\end{proof}
\begin{rem}\label{rem:1}
  We stress the fact that Theorem~\ref{thm:regularityNSE} can be
  stated in a more precise way. In fact, when $\beta \in L_n$, the
  exact condition for the regularity is given
  by~\eqref{eq:regularitiesalpha}$_2$. Hence, another way to state it
  could be the following where we express $\alpha$ as a function of
  $\beta$.
  
  Let $\beta>3$. For every $n\ge 1$, if
\begin{equation*}
    \lambda \in 
    \begin{cases}
      L^{\frac{4(n+1)\beta}{2(n+1)(\beta-3)-\beta}}((0,T;L^\beta(\Omega))
      & \text{if }\beta
      \in R_n
      \\
      L^{\frac{2(n+2)\beta}{2\beta-3}}(0,T;L^\beta(\Omega)) &\text{if }\beta \in L_n,
\end{cases}
\end{equation*}
where $L_n$ and $R_n$ are defined in~\eqref{eq:LnRn},
then the solution $u$ is strong.
\end{rem}

\begin{rem}\label{rem:2}
  From the proof of Theorem~\ref{thm:regularityNSE} we know that at
  every step $n\ge 1$ of the iteration, $\alpha$ and $\beta$ satisfy
  \begin{equation*}
    1-\frac{1}{2(n+1)}\le \frac{2}{\alpha}+\frac{3}{\beta}\le 1-\frac{1}{2(n+2)}.
  \end{equation*}
  Since the left hand side is a strictly increasing function in $n$
  converging to $1$, then we can say that
  \begin{equation*}
    \frac{2}{\alpha}+\frac{3}{\beta} \in \big[3/4,1\big).
  \end{equation*}
  Of course when the previous term is less than $3/4$, regularity
  holds a fortiori.
\end{rem}
The result we proved is the most precise in terms of the exponents,
but is rather difficult for what concern the interpretation, As a direct consequence of Theorem~\ref{thm:regularityNSE}, we have
stated Corollary~\ref{cor:regularityNSE}, which has a clear and simple
statement (even if not being the sharpest result), and  which we can now easily
prove.
\begin{proof}[Proof of Corollary~\ref{cor:regularityNSE}]
  By Remark~\ref{rem:2} we can suppose that
  \begin{equation*}
    \frac{3}{4}\le\frac{2}{\alpha}+\frac{3}{\beta}< 1,
  \end{equation*}
  and, for each couple $(\alpha,\beta)$ satisfying the hypotheses
  of~Corollary\ref{cor:regularityNSE}, there exists a unique
  $\bar{n}\in \mathbb N$ such that
%
  \begin{equation}
    \label{eq:rangebarn}
     1-\frac{1}{2(\overline{n}+1)}\le \frac{2}{\alpha}+\frac{3}{\beta}\le
     1-\frac{1}{2(\overline{n}+2)}. 
   \end{equation}
Now, if we choose
\begin{equation*}
  \beta_0 := \frac{6(\overline{n}+1)}{2\overline{n}-1},
\end{equation*}
it turns out to be the left endpoint of the interval $L_{\overline{n}}$,
and then for 
every $\beta \in (3,\beta_0)$, we have that
\begin{equation*}
     \beta \in \bigcup_{n\ge \overline{n}+1} I_n= (3,\beta_0).
\end{equation*}
In order to conclude we observe that condition~\eqref{eq:rangebarn}
implies the ones given by~\eqref{eq:regularityLnRn} for every
$n\ge \overline{n}+1$, hence concluding the proof.
\end{proof}

\section*{Acknowledgments}
The authors are members of INdAM GNAMPA and they are funded by MIUR
within project PRIN20204NT8W ``Nonlinear evolution PDEs, fluid
dynamics and transport equations: theoretical foundations and
applications'' and MIUR Excellence, Department of Mathematics,
University of Pisa, CUP I57G22000700001.
EC is also funded by  PRIN2022T9K54B ``Classical equations of compressible fluids mechanics: existence and properties of non-classical solutions''.
\section*{Conflicts of interest and data availability statement}
The authors declare that there is no conflict of interest. Data sharing not applicable to this article as no datasets were generated or analyzed during the current study.


\begin{thebibliography}{10}

\bibitem{Abe2022}
K.~Abe.
\newblock Existence of vortex rings in {B}eltrami flows.
\newblock {\em Comm. Math. Phys.}, 391(2):873--899, 2022.

\bibitem{BT2018}
C.~Bardos and E.S. Titi.
\newblock Onsager's conjecture for the incompressible {E}uler equations in
  bounded domains.
\newblock {\em Arch. Ration. Mech. Anal.}, 228(1):197--207, 2018.

\bibitem{BKM1984}
J.T. Beale, T.~Kato, and A.~Majda.
\newblock Remarks on the breakdown of smooth solutions for the $3$-{D} {E}uler
  equations.
\newblock {\em Comm. Math. Phys.}, 94(1):61--66, 1984.


\bibitem{Bei1995a}
H.~Beir\~{a}o~da Veiga, \emph{A new regularity class for the {N}avier-{S}tokes equations in
  {${\bf R}\sp n$}}, Chinese Ann. Math. Ser. B \textbf{16} (1995), no.~4,
  407--412, A Chinese summary appears in Chinese Ann. Math. Ser. A {\bf 16}
  (1995), no.\ 6, 797

\bibitem{BY2019}
H.~Beir\~{a}o~da Veiga and J.~Yang.
\newblock On the energy equality for solutions to {N}ewtonian and
  non-{N}ewtonian fluids.
\newblock {\em Nonlinear Anal.}, 185:388--402, 2019.

\bibitem{Bel1873}
E.~Beltrami.
\newblock Sui principii fondamentali dell'idrodinamica razionale.
\newblock {\em Mem. dell'Accad. Scienze Bologna}, page 394, 1873.

\bibitem{Ber2002a}
L.~C. Berselli, \emph{On a regularity criterion for the solutions to the 3{D}
  {N}avier-{S}tokes equations}, Differential Integral Equations \textbf{15}
  (2002), no.~9, 1129--1137.

\bibitem{Ber2021}
L.~C. Berselli, \emph{Three-dimensional {N}avier-{S}tokes equations for turbulence},
  Mathematics in Science and Engineering, Academic Press, London, [2021]
  \copyright 2021
  

\bibitem{BC2020}
L.~C. Berselli and E.~Chiodaroli, \emph{On the energy equality for the 3{D}
  {N}avier-{S}tokes equations}, Nonlinear Anal. \textbf{192} (2020), 111704,
  24.

\bibitem{BG2024}
L.~C. Berselli and S.~Georgiadis.
\newblock Three results on the energy conservation for the 3{D} {E}uler
  equations.
\newblock {\em NoDEA Nonlinear Differential Equations Appl.}, 31:33, 2024.

\bibitem{BS2024}
L.~C. Berselli and R.~Sannipoli.
  \newblock Velocity-vorticity geometric constraints for the energy
  conservation of 3D ideal incompressible fluids, Technical report, arXiv:2405.08461,  2024.
   	
\bibitem{BDLSV2019}
T.~Buckmaster, C.~de~Lellis, L.~Sz\'{e}kelyhidi, Jr., and V.~Vicol.
\newblock Onsager's conjecture for admissible weak solutions.
\newblock {\em Comm. Pure Appl. Math.}, 72(2):229--274, 2019.



\bibitem{CL2020}

  A.~Cheskidov and X.~Luo, \emph{Energy equality for
    the {N}avier-{S}tokes equations in weak-in-time {O}nsager spaces},
  Nonlinearity \textbf{33} (2020), no.~4, 1388--1403.


\bibitem{CET1994}
P.~Constantin, W.~E, and E.S. Titi.
\newblock Onsager's conjecture on the energy conservation for solutions of
  {E}uler's equation.
\newblock {\em Comm. Math. Phys.}, 165(1):207--209, 1994.

\bibitem{DLS2009}
C.~De~Lellis and Jr.~L. Sz{\'e}kelyhidi.
\newblock The {E}uler equations as a differential inclusion.
\newblock {\em Ann. of Math. (2)}, 170(3):1417--1436, 2009.

\bibitem{Der2020}
L.~De~Rosa.
\newblock On the helicity conservation for the incompressible {E}uler
  equations.
\newblock {\em Proc. Amer. Math. Soc.}, 148(7):2969--2979, 2020.

\bibitem{DE2019}
T.~D. Drivas and G.~L. Eyink, \emph{An {O}nsager singularity theorem for
  {L}eray solutions of incompressible {N}avier-{S}tokes}, Nonlinearity
\textbf{32} (2019), no.~11, 4465--4482.

\bibitem{DN2018}
T.~D. Drivas and H.~Q. Nguyen, \emph{Onsager's conjecture and anomalous
  dissipation on domains with boundary}, SIAM J. Math. Anal. \textbf{50}
  (2018), no.~5, 4785--4811.
  

\bibitem{EP2016}
A.~Enciso and D.~Peralta-Salas.
\newblock Beltrami fields with a nonconstant proportionality factor are rare.
\newblock {\em Arch. Ration. Mech. Anal.}, 220(1):243--260, 2016.

\bibitem{EG2016}
A.~Ern and J.-L. Guermond, \emph{Mollification in strongly {L}ipschitz domains with application
  to continuous and discrete de {R}ham complexes}, Comput. Methods Appl. Math.
  \textbf{16} (2016), no.~1, 51--75.
  

\bibitem{Fri1995}
U.~Frisch.
\newblock {\em Turbulence, The {L}egacy of {A}.{N}.~{K}olmogorov}.
\newblock Cambridge University Press, Cambridge, 1995.

\bibitem{Gal2000a}
G.~P.~Galdi, \emph{An introduction to the {N}avier-{S}tokes initial-boundary value
  problem}, Fundamental directions in mathematical fluid mechanics, Adv. Math.
  Fluid Mech., Birkh\"auser, Basel, 2000, pp.~1--70

\bibitem{Gal2011}
G.~P.~Galdi, \emph{An introduction to the mathematical theory of the
  {N}avier-{S}tokes equations. {S}teady-state problems.}, Springer Monographs
  in Mathematics, Springer-Verlag, New York, 2011.
  
  
\bibitem{GLS2019}
N.~R. Gauger, A.~Linke, and P.~W. Schroeder.
\newblock On high-order pressure-robust space discretisations, their advantages
  for incompressible high {R}eynolds number generalised {B}eltrami flows and
  beyond.
\newblock {\em SMAI J. Comput. Math.}, 5:89--129, 2019.

\bibitem{Ise2018}
P.~Isett.
\newblock A proof of {O}nsager's conjecture.
\newblock {\em Ann. of Math. (2)}, 188(3):871--963, 2018.

\bibitem{Kol1941}
A.N. Kolmogorov, \emph{The local structure of turbulence in incompressible
  viscous fluids for very large {R}eynolds number}, Dokl. Akad. Nauk SSR
  \textbf{30} (1941), 9--13.

\bibitem{Lio1960}
J.-L. Lions, \emph{Sur la r\'egularit\'e et l'unicit\'e des solutions
  turbulentes des \'equations de {N}avier {S}tokes}, Rend. Sem. Mat. Univ.
  Padova \textbf{30} (1960), 16--23
  
\bibitem{LWY2023}
J.~Liu, Y.~Wang, and Y.~Ye.
\newblock Energy conservation of weak solutions for the incompressible {E}uler
  equations via vorticity.
\newblock {\em J. Differential Equations}, 372:254--279, 2023.

\bibitem{NgNgT2019}
Q.-H. Nguyen, P.-T. Nguyen, and B.~Q. Tang.
\newblock Energy equalities for compressible {N}avier-{S}tokes equations.
\newblock {\em Nonlinearity}, 32(11):4206--4231, 2019.

\bibitem{Ons1949}
L.~Onsager.  \newblock Statistical hydrodynamics.  \newblock {\em
    Nuovo Cimento (9)}, 6(Supplemento, 2 (Convegno Internazionale di
  Meccanica Statistica)):279--287, 1949.

\bibitem{Pro1959}
G.~Prodi, \emph{Un teorema di unicit\`a per le equazioni di {N}avier-{S}tokes},
  Ann. Mat. Pura Appl. (4) \textbf{48} (1959), 173--182.

\bibitem{RRS2018}
J.~C. Robinson, J.~L. Rodrigo, and J.~W.~D. Skipper, \emph{Energy conservation
  for the {E}uler equations on {$\mathbb{ T}^2\times\mathbb{R}_+$} for weak
  solutions defined without reference to the pressure}, Asymptot. Anal.
  \textbf{110} (2018), no.~3-4, 185--202.



\bibitem{vWah1992}
W.~von Wahl, \emph{Estimating {$\nabla u$} by {${\rm div}\, u$} and
    {${\rm curl}\, u$}}, Math. Methods Appl. Sci. \textbf{15} (1992),
  no.~2, 123--143.
  
\bibitem{WWWY2023}
Y.~Wang, W.~Wei, G.~Wu, and Y.~Ye.
\newblock On the energy and helicity conservation of the incompressible {E}uler
  equations.
\newblock Technical Report arXiv:2307.08322, 2023.


\end{thebibliography}
\def\ocirc#1{\ifmmode\setbox0=\hbox{$#1$}\dimen0=\ht0 \advance\dimen0
  by1pt\rlap{\hbox to\wd0{\hss\raise\dimen0
  \hbox{\hskip.2em$\scriptscriptstyle\circ$}\hss}}#1\else {\accent"17 #1}\fi}
  \def\polhk#1{\setbox0=\hbox{#1}{\ooalign{\hidewidth
  \lower1.5ex\hbox{`}\hidewidth\crcr\unhbox0}}} \def\cprime{$'$}

\end{document}